\documentclass[11pt]{amsart}
\usepackage{url}
\usepackage{hyperref,amssymb}
\usepackage{appendix}
\hypersetup{colorlinks=true,linkcolor=blue,
urlcolor=cyan,citecolor=magenta} 
\usepackage[hyperpageref]{backref}
\usepackage{frcursive,calrsfs,mathrsfs}
\usepackage{xcolor} 
\newcommand{\Q}{{\mathbb{Q}}}
\newcommand{\Z}{{\mathbb{Z}}}
\newcommand{\F}{{\mathbb{F}}}

\newcommand{\N}{{\mathbb{N}}}
\newcommand{\ds}{\displaystyle}

\newcommand{\ov}{\overline}
\newcommand{\wt}{\widetilde}

\newcommand{\ft}{\footnotesize}
\newcommand{\ns}{\normalsize}
\newcommand{\BN}{{\bf N}}

\newcommand{\BH}{{\bf H}}
\newcommand{\CA}{{\mathcal A}}

\newcommand{\CD}{{\mathcal D}}

\newcommand{\CH}{{\mathcal H}}

\newcommand{\CM}{{\mathcal M}}

\newcommand{\CO}{{\mathcal O}}
\newcommand{\CR}{{\mathcal R}}
\newcommand{\CT}{{\mathcal T}}
\newcommand{\CU}{{\mathcal U}}

\newcommand{\CW}{{\mathcal W}}
\newcommand{\CX}{{\mathcal X}}

\newcommand{\order}{\raise0.8pt \hbox{${\scriptstyle \#}$}}
\newcommand{\lien}{\mathrel{\mkern-4mu}}
\newcommand{\too}{\relbar\lien\rightarrow}

\newcommand{\plus}{\ds\mathop{\raise 0.5pt \hbox{$\bigoplus$}}\limits}
\newcommand{\prd}{\ds\mathop{\hbox{$\prod$}}\limits}
\newcommand{\sm}{\ds\mathop{\raise 1.0pt \hbox{$\sum$}}\limits}
\newcommand{\ffrac}[2]{\hbox{\ft $\displaystyle\frac{#1}{#2}$}}

\newcommand{\Gal}{{\rm Gal}}
\newcommand{\rk}{{\rm rank}}
\newcommand{\tor}{{\rm tor}}

\newcommand{\pr}{{\rm pr}}

\newcommand{\nr}{{\rm nr}}

\newcommand{\lc}{{\rm lc}}

\newcommand{\bp}{{\rm bp}}
\newcommand{\cyc}{{\rm c}}
\newcommand{\acyc}{{\rm ac}}

\newcommand{\Ccl}{c\hskip-0.1pt{\ell}}

\newcommand{\Log}{{\rm Log}}

\newtheorem{theorem}{Theorem}[section]
\newtheorem{lemma}[theorem]{Lemma}

\newtheorem{proposition}[theorem]{Proposition}

\newtheorem{definition}[theorem]{Definition}

\newtheorem{remark}[theorem]{Remark}
\newtheorem{remarks}[theorem]{Remarks}

\numberwithin{equation}{section}
\usepackage[english]{babel}
\title[Test of $p$-rationality of any number field]
{Test of $p$-rationality of a number field \\ pari/gp programs}

\author{Georges Gras}

\address{Villa la Gardette, 4 chemin de Ch\^ateau Gagni\`ere, 
F-38520 Le Bourg d'Oisans}
\email{g.mn.gras@wanadoo.fr}
\urladdr{\url{http://orcid.org/0000-0002-1318-4414}}

\keywords{$p$-rationality, $p$-ramification theory,
class field theory, PARI/GP programs}

\subjclass{11R29, 11R37}

\begin{document}

\date{November 28, 2026}

\begin{abstract}
We give a new computational approach of $p$-rationality of a number field
$k$, when the method using ray class groups is out of reach in a reasonable 
time if $[k : \Q]$ and $p$ are too large. The method is based for a 
part on the utilisation of the $\Log$-function that we have defined in 
the $80$'s. The general {\sc pari/gp} program \ref{P1} makes no assumptions 
on $k$ and $p$. For quadratic fields, we give specific programs for more 
information. Other programs \ref{P5}, \ref{PW}, \ref{P6}, 
calculate the three invariants associated 
to $p$-rationality. Surprisingly, these programs, with many instructions, are 
much faster, and allow to see the behavior of very large $p$.

\bigskip\noindent
{\sc R\'esum\'e.} 
Nous donnons une nouvelle approche num\'erique de la $p$-rationalit\'e d'un 
corps de nombres $k$, lorsque la m\'ethode utilisant des groupes de classes 
de rayon est inaccessible en un temps 
raisonnable si $[k : \Q]$ et $p$ sont tr\`es grands. La m\'ethode
est bas\'ee pour une part sur la fonction $\Log$ que nous avions 
d\'efinie dans les ann\'ees $80$, Le programme {\sc pari/gp} g\'en\'eral 
\ref{P1}  ne fait aucune hypoth\`ese sur $k$ et $p$. Pour les corps 
quadratiques, nous donnons des programmes sp\'ecifiques pour d'avantage 
d'information. D'autres programmes \ref{P5}, \ref{PW}, \ref{P6},
calculent les trois
invariants associ\'es \`a la $p$-rationalit\'e. De fa\c con surprenante, ces 
programmes, avec beaucoup d'instructions, sont beaucoup plus rapides 
et permettent de voir le comportement de tr\`es grands $p$.
\end{abstract}

\maketitle

\tableofcontents

\section{Introduction and Main computational results} 
We intend to furnish various {\sc pari/gp} \cite{Pari2019} programs testing 
the $p$-rationality\,\footnote{\,Under Leopoldt's conjecture, $k$ is $p$-rational 
if the Galois group of the maximal abelian $p$-ramified pro-$p$-extension 
of $k$ is free; in other words its torsion group $\CT_k$ is trivial. 
For more definitions, properties ans story of the 
$p$-rationality as an aspect of the $p$-ramification theory, we refer to 
\cite{Gra1983, Gra1986}, \cite{Jau1986}, \cite{NQD1986}, \cite{Mov1988}, 
\cite{JN1993}, \cite{Jau1998}, \cite{Gra2014}, \cite{Gra2019b}, and their 
bibliographies. We suggest to compute the Jaulent's logarithmic class 
group $\wt \CH_k$ \cite{BJ2016, Jau1994} (instruction ${\sf bnflog(k,p)[1]}$), 
since many similarities exist with $\CT_k$.} 
for any number field $k$ given by means of an irreducible monic 
polynomial $P \in \Z[x]$; the first one (Program \ref{P1}) is new, 
completely general and runs much faster. 
Except for few ``bad primes'', the program does not use ray class groups.
For comparison, we give again the classical Program \ref{P2} only
using ray class groups of suitable modulus (adapted version of \cite[Program II, 
p. 37]{Gra2019a}); its execution time is exponential with respect to~$p$,
and it does not give the details in the formula $\order \CT_k = \order \CH_k^\infty 
\cdot \CR_k \cdot \order \CW_k^\bp$ (cf. \S\,\ref{pramif}).
The two others Programs \ref{P3}, \ref{P4} are
specified for quadratic fields and do not use ray class groups.

\smallskip
In general, we provide only simple versions of the instructions of the 
{\sc pari/gp} language, so that any reader (whether or not they use 
{\sc pari/gp}) can analyze the programs and make use of them.
The programs may be copied and pasted from any ${\sf tex}$ file (easier) 
or ${\sf pdf}$ file (copy separately onto each page); if necessary, the reader 
can download the {\sc par/gp}-package 
from \url{https://pari.math.u-bordeaux.fr/gpwasm.html}.\,\footnote{\,All our
programs are suitable for the use of the ``Run PARI/GP in your browser'', 
but we are warned that ``speed is approximately 15\% of native speed''.}

\smallskip
I apologize if certain programs appear to be causing issues, but all 
them have been as confirmed as possible. I would be 
grateful to anyone who could provide me with information on this subject. 
It is always possible to check using the very safe Program \ref{P2} involving
ray class groups, provided that $p$ is small; but our purpose is to 
analyse, for some generalized conjectures (e.g. \cite{Gra2016}), the huge 
$p$ giving a non $p$-rationality (if any). Moreover, since a $p$-rational
field has the simplest arithmetic regarding its $p$-invariants, it is necessary 
to be able to test any ``random'' $k$ contrary to an abundant literature
giving criteria for ``too classical'' number fields.

\smallskip
In the next section, we give Program \ref{P5} for the test of triviality of the 
normalized $p$-adic regulator $\CR_k$.\footnote{\,Recall that $\CR_k(p)$
is also, by abuse of notation, the order in $p^\N$ of a $p$-adic regulator 
(= matrix determinant) and that $\order \CR_k(p) = 1$, as $\Z_p$-module, 
means $\CR_k(p)$ invertible in $\Z_p$ as determinant.} Then we give 
Programs \ref{PW1}, \ref{PW2}, \ref{PW3}, testing, with various 
approaches, the condition $\CW_k^\bp := \bigoplus_{v \mid p} 
\mu_p(k_v)\big/\mu_p(k) = 1$, in terms of local 
and global $p$-groups of roots of unity. Finally, Program \ref{P6} computes
the part of the $p$-Hilbert class field of $k$ contained in the 
compositum of the $\Z_p$-extensions of $k$, and deduces the tricky
invariant $\CH_k^\infty$, kernel of the $\Log$ map (see \S\,\ref{Logp}).

\subsection{General program for any number field and any prime number}
Refer to the Section \ref{pramif} for some class field theory justifications and 
the definition of the finite $\Z_p$-modules $\CH_k^\infty$, $\CR_k$, $\CW_k^\bp$ 
(in ${\sf Hinfty}$, ${\sf MLogE}$, ${\sf W}$, respectively in 
the programs)\,\footnote{\,In 
the paper, when we evoke programmatic aspects, we will use some verbatim 
notations that appear in the {\sf pari/gp} programs, as for instance ${\sf k}$, 
${\sf P}$, ${\sf d}$, ${\sf r1}$, ${\sf r2}$, ${\sf p}$, ${\sf Dk}$, ${\sf Ek}$, 
${\sf Hk}$, ${\sf Wk}$, ${\sf Rk}$, etc.}, 
being such that the $p$-rationality of $k$ is equivalent to 
$\CH_k^\infty = \CR_k = \CW_k^\bp = 1$.
Program \ref{P1} is not based on computations of ray class groups,
except for few bad primes that otherwise would need very complicate 
sub-algorithms; for the non-bad primes, computations are much 
easier and faster.

\subsubsection{Bad primes associated to the field $k$}\label{pBad}

See Definition \ref{bad} for more details on the notion of bad primes for 
which the programs use classical way. These primes are few in number 
and generally small. 
The following program gives the prime numbers $p$ dividing the class 
number ${\sf hk=k.no}$, the discriminant ${\sf Dk=k.disc}$, the denominators 
${\sf denB}$ of the integer basis elements ${\sf k.zk[i]}$ of $k$, and the 
primes less than $d+1$, where $d = [k : \Q]$:

\smallskip
\ft\begin{verbatim}
allocatemem(800000000)
{P=x^29+7*x^3-14*x-7;d=poldegree(P);if(polisirreducible(P)==0,break);
print("P=",P);k=bnfinit(P,1);print("Hk=",k.cyc);Pbad=List;
forprime(p=2,d+1,listput(Pbad,p));
denB=1;for(j=1,d,Bzk=lift(k.zk[j]);for(i=0,d-1,
denB=lcm(denominator(polcoeff(Bzk,i)),denB)));
\\PRIME DIVISORS OF hk,Dk,denB:
Fbad=factor(k.no*abs(k.disc)*denB);
dF=matsize(Fbad)[1];for(i=1,dF,p=component(Fbad,1)[i];
if(p>d+1,listput(Pbad,p)));print("Pbad=",Pbad)}
time = 14,763 ms.

P=x^9+7*x^3-14*x-7   Hk=[]
Pbad=List([2,3,5,7,251,4649])
P=x^19+7*x^3-14*x-7  Hk=[]
Pbad=List([2,3,5,7,11,13,17,19,1158587,290260991,335367508213])
P=x^29+7*x^3-14*x-7  Hk=[]
Pbad=List([2,3,5,7,11,13,17,19,23,29,67,251,265403,8929674616619,55306100356254142133177])
\end{verbatim}\ns

\begin{remarks}
(i) One verifies that primes dividing denominators of units, divide that 
of integer basis elements.

\smallskip
(ii) The huge bad primes will probably not be considered by the selected
intervals in the programs; they are due in general to ramified primes as for 
the above polynomial where the field discriminant is
\ft${\sf  Dk = 7^{28} \cdot 251 \cdot 265403 \cdot 8929674616619 \cdot 55306100356254142133177}$.\ns

\smallskip
For such huge primes, one checks easily that $\CH_k^\infty = \CW_k^\bp = 1$, 
but $\CR_k \ne 1$ remains a possible obstruction (Borel--Cantelli probability 
very close to $0$). If the prime $p$ is not bad but too large for a program, 
one may use the specific Program \ref{P5bis} testing only the regulator.

\smallskip
(iii) The primes less than or equal to $d+1$ are considered as bad primes 
since they may give $\CW_k^\bp \ne 1$ because $\mu_p \ne 1$ (in $k$ 
or $k_v$) implies $[\Q_p(\mu_p) : \Q_p] = p-1 \leq d$.
\end{remarks}

\subsubsection{Program P1}\label{P1}
The reader has only to choose the polynomial ${\sf P}$ and the interval
${\sf [bp, Bp]}$ at will. The upper bound ${\sf Maxbad}$ avoid the huge
bad primes and may be modified. Since the non-bad primes are odd and 
unramified by definition, the maximal summation index giving
$\log$ series modulo ${\sf p^N}$, ${\sf N \in \{2, 3\}}$ may be 
${\sf hatt=N}$ (see Proposition \ref{prop1}).

\smallskip
Unless otherwise stated, the polynomials, in our examples, have been 
chosen at random, but it is then very difficult to find non $p$-rational 
examples when $\mu_p(k) = 1$:

\smallskip
\ft\begin{verbatim}
{P=x^9+7*x^3-14*x-7;d=poldegree(P);bp=2;Bp=10^6;Maxbad=5*10^3;
if(polisirreducible(P)==0,break);print("P=",P);k=bnfinit(P,1);
print("class group=",k.cyc);kP=nfinit(P);rk=kP[2][2]+1;
print("number of independent Z_p-extensions=",rk);
E=k.fu;rE=kP[2][1]+kP[2][2]-1;print("Z-rank of the unit group=",rE);
\\BAD PRIMES:
Pbad=List;
\\SMALL BAD PRIMES:
Minbad=d+1;dbad=0;forprime(p=bp,Minbad,listput(Pbad,p);dbad=dbad+1);
\\DENOMINATORS OF THE INTEGER BASIS:
denB=1;for(j=1,d,Bzk=lift(k.zk[j]);for(i=0,d-1,
denB=lcm(denominator(polcoeff(Bzk,i)),denB)));
\\PRIME DIVISORS OF hk,Dk,denB:
Fbad=factor(k.no*abs(k.disc)*denB);dF=matsize(Fbad)[1];gbad=0;
for(i=1,dF,p=component(Fbad,1)[i];if(p>Minbad,listput(Pbad,p);
gbad=gbad+1));print("Pbad=",Pbad);Nbad=dbad+gbad;
\\PRODUCT OF ALL THE BAD PRIMES:
Prodbad=1;for(q=1,Nbad,Prodbad=Prodbad*Pbad[q]);
\\TEST OF p-RATIONALITY FOR THE BAD PRIMES < Maxbad:
for(j=1,Nbad,p=Pbad[j];if(p>Maxbad,next);
\\USE OF A RAY CLASS GROUP FOR Pbad[j]:
nup=2;if(p==2,nup=3);Tk=List;L=List;kpnup=bnrinit(k,p^nup);
Hpnup=kpnup.cyc;e=component(matsize(Hpnup),2);R=0;
for(j=1,e,c=Hpnup[e-j+1];if(Mod(c,p)==0,R=R+1;
listput(L,p^valuation(c,p))));
\\WRITTINGS: p-RATIONALITY/NON-p-RATIONALITY:
\\SUPPRESS \\ BELOW TO WRITE ALL p-RATIONAL FIELDS:
\\if(R==rk,print("k is ",p,"-rational"));
if(R>rk,for(i=1,R-rk,c=L[i];listput(Tk,c));
print("p=",p," rank of the ",p,"-torsion group Tk=",R-rk);
print("k is not ",p,"-rational");
print("**********************")));
\\TEST FOR THE NON-BAD PRIMES < Bp:
\\DECOMPOSITION OF p IN k, COMPUTATION OF M:
forprime(p=bp,Bp,if(Mod(Prodbad,p)==0,next);
F=idealfactor(k,p);g=matsize(F)[1];M=1;for(i=1,g,
Primei=component(F,1)[i];fi=Primei[4];M=lcm(M,p^fi-1));
\\TEST OF INVERSIBILITY OF THE REGULATOR:
N=2;mod=p^N;hatt=N;LogE=List;
\\REDUCTION OF epsilonM MODULO mod:
for(ell=1,rE,epsilonM=Mod(Mod(lift(E[ell]),mod),P)^M;
\\COMPUTATION OF Logarithm OF epsilonM:
pi=Mod(lift(lift(epsilonM)),P)-1;Log=pi;t=1;pit=pi;
while(t<=hatt,t=t+1;pit=pi*pit;Log=Log+(-1)^(t+1)*pit/t;
\\REDUCTION OF Log MODULO mod:
Log=Mod(Mod(lift(Log),mod),P));
Log=Mod(lift(lift(Log)),P);listput(LogE,Log));
\\MATRIX MLogE ((rE,d)-matrix) GIVING Rk MODULO p:
ME=List;for(ell=1,rE,Lell=LogE[ell];Lell=lift(Lell)/p;
Lc=List;LellZ=nfalgtobasis(k,Mod(Lell,P));
for(j=1,d,coeff=Mod(LellZ[j],p);listput(Lc,coeff));listput(ME,Lc));
MLogE=matrix(rE,d,u,v,ME[u][v]);
\\COMPUTATION OF THE Z/pZ-rank OF MLogE:
RLog=matrank(MLogE);
\\WRITTINGS: p-RATIONALITY/NON-p-RATIONALITY:
\\if(RLog==rE,print("k is ",p,"-rational"));
if(RLog<rE,print("p=",p," rank(MLogE)=",RLog);
print("k is not ",p,"-rational");
print("**********************")))}

P=x^9+7*x^3-14*x-7
class group=[]
number of independent Z_p-extensions=4
Z-rank of the unit group=5
Pbad=List([2,3,5,7,251,4649])
**********************
p=2 rank(MLogE)=1
k is not 2-rational
**********************
**********************
p=61 rank(MLogE)=4
k is not 61-rational
**********************
**********************
p=5419 rank(MLogE)=4
k is not 5419-rational
**********************
**********************
p=12703 rank(MLogE)=4
k is not 12703-rational
**********************
for p<10^4: time = 720 ms.
for p<10^5: time = 8,345 ms.
for p<10^6: time = 1min, 34,310 ms.
for p<10^7: time = 15min, 44,164 ms.
\end{verbatim}\ns

There is no other solutions up to $10^9$.

\smallskip
This example has the curious property that the integer basis ${\sf k.zk}$
is given by $9$ units, including the set of the $5$ fundamental units ${\sf E[i]}$:

\smallskip
\ft\begin{verbatim}
nfalgtobasis(k,E[1])     =[0,0,0,0,-1,0,0,0,0]
nfalgtobasis(k,E[2])     =[0,0,0,0,0,0,-1,0,0]
nfalgtobasis(k,E[3])     =[0,0,0,0,0,1,0,0,0]
nfalgtobasis(k,E[4])     =[0,0,1,1,-1,-1,0,0,0]
nfalgtobasis(k,E[5])     =[0,0,0,1,-1,0,0,0,0]
nfalgtobasis(k,E[1]*E[2])=[0,-1,-1,0,0,0,0,0,0]
nfalgtobasis(k,E[2]*E[3])=[1,0,0,0,0,0,1,0,0]
nfalgtobasis(k,E[3]*E[4])=[-1,-1,0,0,0,0,0,-1,0]
nfalgtobasis(k,E[4]*E[5])=[-1,0,0,-1,1,0,0,-1,1]
\end{verbatim}\ns

\subsection{Classical program testing \texorpdfstring{$p$}{Lg}-rationality 
with ray class groups} 
This program is using the structure of a suitable ray class group:

\subsubsection{Program P2}\label{P2}
The values ${\sf nup = 2}$ for $p>2$, ${\sf nup = 3}$ for $p = 2$, are 
sufficient for computing the $p$-rank of the $p$-torsion group $\CT_k$; 
otherwise, take a larger value of ${\sf nup}$ to get its complete 
structure (knowing the $p$-rank of ${\sf Tk}$, it suffices to stabilize
the corresponding ${\sf R-rk}$ components given by ${\sf kpnup.cyc}$,
the ${\sf rk}$ first ones, being due to the $\Z_p$-extensions, and tend to 
infinity with ${\sf nup}$):

\smallskip
\ft\begin{verbatim}
allocatemem(1000000000)
{\\P=x^4+8*x^2+324;
P=x^9+7*x^3-14*x-7;
if(polisirreducible(P)==0,break);bp=2;Bp=10^4;
print("P=",P);rk=nfinit(P)[2][2]+1;
print("number of independent Z_p-extensions: ",rk);
k=bnfinit(P,1);print("class group=",k.cyc);
\\PRIMES p AND CHOICE OF nup:
forprime(p=bp,Bp,Tk=List;L=List;
nup=2;if(p==2,nup=3);
\\COMPUTATION OF THE p-RANK OF THE TORSION GROUP Tk:
kpnup=bnrinit(k,p^nup);Hpnup=kpnup.cyc;
e=component(matsize(Hpnup),2);R=0;for(j=1,e,c=Hpnup[e-j+1];
if(Mod(c,p)==0,R=R+1;listput(L,p^valuation(c,p))));
if(R>rk,for(j=1,R-rk,c=L[j];listput(Tk,c));print();
\\STRUCTURE OF Tk[p]:
print("rk(Tk)=",R-rk," k is not ",p,"-rational"," Tk[p]=",Tk)))}

P=x^9+7*x^3-14*x-7
number of independent Z_p-extensions: 4
class group=[]
rk(Tk)=1 k is not 2-rational      Tk[p]=List([2])
rk(Tk)=1 k is not 61-rational     Tk[p]=List([61])
rk(Tk)=1 k is not 5419-rational   Tk[p]=List([5419])
time = 201 ms for p<10^2.
time = 17,130 ms for p<10^3.
time = 1h, 35min, 7,927 ms for p<10^4.
\end{verbatim}\ns

The program is very slow for large primes and must be used only for 
small ones; moreover, one must increase the memory (e.g., ${\sf 
allocatemem(1000000000)}$). It is difficult to compute beyond $10^3$;
on the contrary, Program \ref{P1} takes ${\sf time = 15min, 42,660 ms}$ 
up to $10^7$ !

\subsection{Specific program for imaginary quadratic fields}
The program is without use of a ray class group. The method consists of 
calculating separately $\CW_k^\bp(p)$ (immediate tests limited to $2$ and 
$3$ since $[\Q_p(\mu_p) : \Q_p] = p-1$ must divide local degrees $1$ or $2$), 
and $\CH_k^\infty$ (more difficult 
since we must compute the part of the Hilbert class field $H_k^\nr$ contained 
in the compositum $\wt k$ of the $\Z_p$-extensions of $k$, which uses the 
$p$-adic $\Log$-function, whence the ``Fermat quotients'' of the pseudo-units 
given by ideal classes); in fact, for $p>2$, since complex conjugation acts on
the $p$-invariants of $k$, $H_k^\nr \cap \wt k \subset k^\acyc$, 
the anti-cyclotomic $\Z_p$-extension.
Note that in the programs for quadratic fields, if ${\sf pi}$ is any integer of $k$ 
of non-zero ${\mathfrak p}$-valuations for all $p$-places ${\mathfrak p}$, then 
${\sf log(1+pi)}$ is a $p$-integer for all $p \geq 2$, which is not the case in some 
framework of ramification of $p$ (this explains that ramified primes are ``bad'' ones). 
If $p \geq 2$ is unramified, then ${\sf log(1+pi)}$ is a multiple of $p$ (see
Proposition \ref{prop1}).

\subsubsection{Program P3 for imaginary quadratic fields and $p > 2$}\label{P3}

The case $p=2$ will be treated separately (Program \ref{P3bis} below).
Once the test on ${\sf W}$ is done and since there is no regulator, only
fields with a non-trivial $p$-class group are tested. Since $p>2$ and
$k$ quadratic, the ${\mathfrak p}$-valuation function of terms of the 
form $\pi^{p^x}/p^x$ is increasing from $\hat x_{\mathfrak p}<1$ 
(Proposition \ref{prop1}).
The reader has to indicate the intervals ${\sf [minm : maxm]}$,
$m \geq 2$, and ${\sf [minp : maxp]}$, $p \geq 3$:

\smallskip
\ft\begin{verbatim}
allocatemem(800000000)
\p 600
{minm=2;maxm=10^4;minp=3;maxp=10^4;
for(m=minm,maxm,if(core(m)!=m,next);P=x^2+m;k=bnfinit(P);
\\COMPUTATION OF W FOR p=3:
W=1;if(-m!=-3 & Mod(-m,9)==-3,W=3;print("-m=",-m,
" k non ",3,"-rational"," W=",3));
\\TESTS ON THE p-CLASS GROUP RANK:
Hk=k.clgp;Nk=Hk[2];A=Hk[3];dA=matsize(A)[2];
forprime(p=3,Hk[1],valhkp=valuation(Hk[1],p);if(valhkp==0,next);
dp=0;for(i=1,dA,if(Mod(Nk[i],p)==0,dp=dp+1));if(dp>1,print("-m=",-m,
" k non ",p,"-rational, rankHk=",dp," W=",W," Hk=",Nk);next);
\\THE p-CLASS GROUP IS NOW p-CYCLIC (A0,N0):
\\FINDING OF THE UNIQUE COMPONENT OF THE p-CLASS GROUP:
for(i=1,dA,if(Mod(Nk[i],p)==0,N0=Nk[i];A0=A[i];break));
n0=valuation(N0,p);N=n0+2;mod=p^N;
\\DECOMPOSITION OF p:
F=idealfactor(k,p);g=matsize(F)[1];Prime=component(F,1)[1];
\\PRIME-TO-p IDEAL GENERATOR B0 OF THE p-CLASS OF A0:
B0=A0;Norm0=idealnorm(k,B0);w0=valuation(Norm0,p);
if(w0!=0,
\\RAMIFIED CASE:
if(Mod(m,p)==0,z=idealappr(k,B0);B0=idealdiv(k,z,B0));
\\UNRAMIFIED CASE:
if(Mod(m,p)!=0,Bp=idealmul(k,B0,p);z=idealappr(k,Bp);
B0=idealdiv(k,z,Bp)));
\\COMPUTATION of THE VALUATION OF Log(U)^-: 
if(p>3,
if(Mod(-m,p)!=0,vLogU=1);if(Mod(-m,p)==0,vLogU=0));
if(p==3,
if(Mod(-m,3)!=0,vLogU=1);if(Mod(-m,3)==0 & Mod(-m,9)!=-3,vLogU=0);
if(Mod(-m,3)==0 & Mod(-m,9)==-3,vLogU=1));
\\DECOMPOSITION OF (p)--COMPUTATION OF M (p^f-1):
rho=Prime[3];f=Prime[4];M=p^f-1;
\\COMPUTATION OF beta0 AND alpha0M:
C0=idealpow(k,B0,N0);beta0=bnfisprincipal(k,C0)[2];
alpha0=Mod(beta0[1]*k.zk[1]+beta0[2]*k.zk[2],P);
\\REDUCTION OF alpha0^M MODULO mod:
alpha0M=Mod(Mod(lift(alpha0),mod),P)^M;
pi=Mod(lift(lift(alpha0M))-1,P);if(pi==0,pi=Mod(mod,P));
\\COMPUTATION OF hatt FOR THE LOG SERIES:
Delta=idealval(k,pi,Prime);
n=1;while(p^n*Delta-(N+n)*rho<0,n=n+1);
hatn=n;if(hatn==1,hatt=N);if(hatn>1,hatt=p^hatn;
while(hatt*Delta-(valuation(hattj,p)+N)*rho>0,hatt=hatt-1));
\\COMPUTATION OF Log0=Log(alpha0M):
t=1;pit=pi;Log0=pi;while(t<=hatt,t=t+1;pit=pi*pit;
Log0=Log0+(-1)^(t+1)*pit/t;
\\REDUCTION OF Log0 MODULO mod:
Log0=Mod(Mod(lift(Log0),mod),P));
a=polcoeff(lift(Log0),1);if(a==0,vala=N);if(a!=0,vala=valuation(a,p));
\\VALUATION OF LogB0 AND COMPUTATION OF Index::
vlogB0=vala-n0;Minval=min(vLogU,vlogB0);
Index=vLogU-Minval;valhkinfty=valhkp-Index;if(valhkinfty>0,
print("-m=",-m," k non ",p,"-rational"," W=",W," Index=",p^Index,
" hkinfty=",p^valhkinfty," Hk=",Nk))))}
time = 26,767 ms.
\end{verbatim}\ns 

There are many non-$p$-rational fields fields for $m \in [2, 10^4]$ and
$p \in [3,10^4]$ but essentially when $\CW_k^\bp \ne 1$ and when 
$\CH_k^\infty \ne 1$, $\rk_p(\CH_k) > 2$, for $p = 3, 5, 7, 11, 13, 17, 19, 23$.

\smallskip
For $p \in [29, 10^6]$, one finds:

\smallskip
\ft\begin{verbatim}
-m=-4106  k non 29-rational Index=1  hkinfty=29  Hk=[58]
-m=-4682  k non 31-rational Index=1  hkinfty=31  Hk=[62]
-m=-5333  k non 31-rational Index=1  hkinfty=31  Hk=[62]
-m=-8597  k non 31-rational Index=1  hkinfty=31  Hk=[62]
-m=-9203  k non 31-rational Index=1  hkinfty=31  Hk=[31]
\end{verbatim}\ns

For $p = 3$ and a cyclic $3$-class group such that
$H_k^\nr \subset k^\acyc$, we obtain, up to $10^4$:

\smallskip
\ft\begin{verbatim}
-m=-1879 k non 3-rational   Index=9   hkinfty=3 Hk=[27]
-m=-2759 k non 3-rational   Index=9   hkinfty=3 Hk=[54]
-m=-3299 k non 3-rational   Index=9   hkinfty=3 Hk=[9,3]
-m=-3671 k non 3-rational   Index=27  hkinfty=3 Hk=[81]
-m=-3943 k non 3-rational   Index=9   hkinfty=3 Hk=[27]
-m=-4119 k non 3-rational   Index=9   hkinfty=3 Hk=[54]
-m=-4397 k non 3-rational   Index=9   hkinfty=3 Hk=[54]
-m=-4541 k non 3-rational   Index=9   hkinfty=3 Hk=[54,2]
-m=-5703 k non 3-rational   Index=9   hkinfty=3 Hk=[18,3]
-m=-7661 k non 3-rational   Index=9   hkinfty=3 Hk=[54,2]
-m=-7703 k non 3-rational   Index=27  hkinfty=3 Hk=[81]
-m=-8426 k non 3-rational   Index=9   hkinfty=3 Hk=[54,2]
-m=-8506 k non 3-rational   Index=9   hkinfty=3 Hk=[54]
-m=-8807 k non 3-rational   Index=27  hkinfty=3 Hk=[81]
-m=-8867 k non 3-rational   Index=9   hkinfty=3 Hk=[27]
-m=-8951 k non 3-rational   Index=9   hkinfty=3 Hk=[135]
-m=-9031 k non 3-rational   Index=9   hkinfty=3 Hk=[54]
-m=-9074 k non 3-rational   Index=9   hkinfty=3 Hk=[54,2]
-m=-9574 k non 3-rational   Index=9   hkinfty=3 Hk=[18,3]
-m=-9663 k non 3-rational   Index=9   hkinfty=3 Hk=[54]
-m=-9749 k non 3-rational   Index=2   hkinfty=3 Hk=[162]
-m=-9831 k non 3-rational   Index=9   hkinfty=3 Hk=[54,2]
\end{verbatim}\ns

\subsubsection{Complementary Program P3 for $p = 2$}\label{P3bis}
The case $p=2$ is particular because of the computation of $\CH_k^\infty$;
so, in that case the use of the classical way is preferable since the
modulus of the ray class group is only $8$. 
Before we compute $\order \CW_k^\bp(2)$:

\smallskip
\ft\begin{verbatim}
{minm=2;maxm=10^4;for(m=minm,maxm,if(core(m)!=m,next);
P=x^2+m;k=bnfinit(P);rk=2;
\\COMPUTATION OF W:
if(-m!=-1 & (Mod(-m,8)==1 || Mod(-m,8)==-1),print("-m=",-m,
" k is not ",2,"-rational"," W=",2));
\\COMPUTATION OF THE 2-RANK OF Tk:
kpnup=bnrinit(k,8);Hpnup=kpnup.cyc;e=component(matsize(Hpnup),2);
R=0;for(j=1,e,c=Hpnup[e-j+1];if(Mod(c,2)==0,R=R+1));if(R>rk,print("-m=",-m,
" k is not ",2,"-rational"," rk(Tk)=",R-rk)))}
time = 6,209 ms.

-m=-7    k is not 2-rational W=2
-m=-14   k is not 2-rational rk(Tk)=1
-m=-15   k is not 2-rational W=2
-m=-17   k is not 2-rational W=2
-m=-21   k is not 2-rational rk(Tk)=1
(...)
-m=-1785 k is not 2-rational W=2
-m=-1785 k is not 2-rational rk(Tk)=3
(...)
-m=-9993 k is not 2-rational W=2
-m=-9994 k is not 2-rational rk(Tk)=1
-m=-9995 k is not 2-rational rk(Tk)=1
-m=-9997 k is not 2-rational rk(Tk)=1
-m=-9998 k is not 2-rational rk(Tk)=1
\end{verbatim}\ns

\subsubsection{Classical program, for imaginary quadratic fields, 
using ray class groups} \label{P3ter}
The classical program using ray class groups takes much time; it gives the 
same results. An interest is that increasing ${\sf nup}$ (giving a modulus 
${\sf p^{nup}}$ large enough), one obtains the group structure of ${\sf Tk}$
(but without giving the orders of the three factors):

\smallskip
\ft\begin{verbatim}
{bp=2;Bp=10^4;for(m=2,10^4,if(core(m)!=m,next);P=x^2+m;
k=bnfinit(P);forprime(p=bp,Bp,nup=2;if(p==2,nup=3);
kpnup=bnrinit(k,p^nup);Hpnup=kpnup.cyc;e=component(matsize(Hpnup),2);
R=0;for(j=1,e,c=Hpnup[e-j+1];if(Mod(c,p)==0,R=R+1));if(R>2,
print("-m=",-m," Hk=",k.cyc," k is not ",p,"-rational"," rk(Tk)=",R-2))))}
time = 28,361 ms for p<10^2.
time = 25min, 48,895 ms  for p<10^4.
\end{verbatim}\ns

As for a random polynomial the use of ray class groups is not 
satisfactory for imaginary quadratic fields, even if there is no regulator.

\subsection{Specific programs for real quadratic fields}
There is some similarities with the imaginary case, the pseudo-units
$\alpha_i$ being replaced by the fundamental unit $\varepsilon$; in the 
real case, $\CH_k^\infty = \CH_k$, except for $p=2$, when 
$m \equiv 2 \pmod 8$, $m \ne 2$, in which case, $\order \CH_k^\infty 
= \frac{1}{2}\order\CH_k$:

\subsubsection{Program P4 for real quadratic fields}\label{P4}
The program, without using ray class groups, is valid for any 
$p \geq 2$ and test successively $\CH_k^\infty$, $\CW_k^\bp$ 
for $p=2, 3$, then $\CH_k^\infty$.

\smallskip
\ft\begin{verbatim}
{minm=2;maxm=10^4;minp=2;maxp=10^4;
for(m=minm,maxm,if(core(m)!=m,next);P=x^2-m;k=bnfinit(P,1);
forprime(p=minp,maxp,N=2;if(p==2,N=3);mod=p^N;
\\COMPUTATION OF Hinfty:
Hinfty=p^(valuation(k.no,p));if(p==2 & m!=2 & Mod(m,8)==2,
Hinfty=Hinfty/2);
\\COMPUTATION OF W FOR p=2,3:
W=1;
if(p==2 & (Mod(m,8)==1 || Mod(m,8)==-1),W=2);
if(p==3 & Mod(m,9)==-3,W=3);
if(p==2 & (W==2 || Mod(Hinfty,2)==0),print("m=",m,
" k non ",p,"-rational"," W=",W," Hinfty=",Hinfty);next);
if(p==3 & (W==3 || Mod(Hinfty,3)==0),print("m=",m,
" k non ",p,"-rational"," W=",W," Hinfty=",Hinfty);next);
if(p>3 & Mod(Hinfty,p)==0,print("m=",m,
" k non ",p,"-rational"," W=",W," Hinfty=",Hinfty);next);
\\COMPUTATION OF THE p-ADIC REGULATOR:
F=idealfactor(k,p);g=matsize(F)[1];
Prime=component(F,1)[1];rho=Prime[3];f=Prime[4];
\\COMPUTATION OF M (p-1, or p+1 using Norm(epsilon^2)=1):
if(f==1,M=p-1);if(f==2,M=p+1);
\\UNIT OF NORM 1 WITH p-INTEGER COEFFICIENTS:
epsilon=k.fu[1];if(p!=2,alpha=epsilon^2);if(p==2,alpha=epsilon^(3*f));
\\COMPUTATION OF alpha^M MODULO mod:
alphaM=Mod(Mod(lift(alpha),mod),P)^M;
\\COMPUTATION OF hatt AND Log(alphaM):
pi=Mod(lift(lift(alphaM))-1,P);if(pi==0,pi=Mod(mod,P));
if(p==2,hatt=4);if(p>2,hatt=2);Log=pi;t=1;pit=pi;
while(t<=hatt,t=t+1;pit=pi*pit;Log=Log+(-1)^(t+1)*pit/t;
\\REDUCTION OF Log MODULO mod:
Log=Mod(Mod(lift(Log),mod),P));
Log=lift(lift(Log));if(Log==0,Log=Mod(mod,P));
\\TEST OF NON-p-RATIONALITY FOR p>2:
if(p!=2,if(Mod(m,p)==0 & valuation(Log,p)>0,
print("m=",m,"  k non ",p,"-rational",
" W=",W," Hinfty=",Hinfty));
if(Mod(m,p)!=0 & valuation(Log,p)>1,
print("m=",m,"  k non ",p,"-rational",
" W=",W," Hinfty=",Hinfty)));
\\TEST OF NON-2-RATIONALITY, COMPUTATION OF rd:
if(p==2,rd=1;if(Mod(m,8)==2,rd=0);if(Mod(m,8)==5,rd=2);
if(valuation(Log,2)>rd,print("m=",m,
" k non ",p,"-rational"," W=",W," Hinfty=",Hinfty)))))}
time = 2min, 42,705 ms.
\end{verbatim}\ns

\smallskip
It follows from the particular cases of definition of $\alpha$ that the 
non-$p$-rationality needs to include at the end of the program, 
the following instructions, noting that for the unit $\alpha$, 
$\BN_{k/\Q}(\alpha) = 1$ implies $\log(\alpha^M) = a \cdot x$ 
(for {\sc pari/gp}, $x$ is the polynomial indeterminate whose 
arithmetical value is $\sqrt m$):

\smallskip
(i) The case $p \ne 2$ is related to the ramified case when 
$\log(\alpha^M) = a \cdot \pi \sim a \cdot \sqrt m$
is of non-zero valuation, which means $v_p(a) \geq 1$:

\smallskip
\ft\begin{verbatim}
if(p!=2,if(Mod(m,p)==0 & valuation(Log,p)>0,print("m=",m,
" k non-p-rational for p=",p);next);
if(Mod(m,p)!=0 & valuation(Log,p)>1,print("m=",m,
" k non-p-rational for p=",p);next))
\end{verbatim}\ns

\smallskip
(ii) The case $p = 2$ takes into account the interpretation of the 
valuation of $\log(\alpha^M) = a \cdot x$; the ramified case 
$m \equiv 2 \pmod 8$ is similar as the previous one. The case
$m \equiv 5 \pmod 8$ corresponds to $2$ inert and 
$\order\CR_k \sim \frac{1}{2} \log(\varepsilon)$, but in that case
$\alpha = \varepsilon^6$ and a simple calculation shows that 
$\pi \in 4 \CO_k$ and $\log(\alpha^M) \in 4 \CO_k$, whence a 
non-trivial regulator if and only if $v_2(\log(\alpha^M))> 2$. 
Otherwise, the condition is $v_2(\log(\alpha^M))> 0$ 
since $\pi \in 2  x \,\CO_k$, hence $v_2(a) >1$:

\smallskip
\ft\begin{verbatim}
if(p==2,rd=1;if(Mod(m,8)==2,rd=0);if(Mod(m,8)==5,rd=2);
if(valuation(Log,2)>rd,print("m=",m," k non ",p,"-rational);next))
\end{verbatim}\ns

\smallskip
These tests on ${\sf valuation(Log,p)}$ imply that all the reductions 
modulo $p^N$ are with $N = 2$ for $p \ne 2$ and $N=3$ for $p = 2$.
 
\smallskip
One obtains $17253$ non-$p$-rational examples. Contrary
to the imaginary case, the non $p$-rationality can occur
for many primes $p$; we give the following excerpt
when the non-$p$-rationality is due to $\CR_k$:

\smallskip
\ft\begin{verbatim}
m=37    k non-p-rational for p=7    W=1 Hinfty=1
m=37    k non-p-rational for p=89   W=1 Hinfty=1
m=37    k non-p-rational for p=257  W=1 Hinfty=1
m=37    k non-p-rational for p=631  W=1 Hinfty=1
m=42    k non-p-rational for p=2    W=1 Hinfty=1
m=42    k non-p-rational for p=3    W=3 Hinfty=1
m=42    k non-p-rational for p=5    W=1 Hinfty=1
m=42    k non-p-rational for p=43   W=1 Hinfty=1
m=42    k non-p-rational for p=71   W=1 Hinfty=1
m=109   k non-p-rational for p=3    W=1 Hinfty=1
m=109   k non-p-rational for p=5    W=1 Hinfty=1
m=109   k non-p-rational for p=11   W=1 Hinfty=1
m=109   k non-p-rational for p=19   W=1 Hinfty=1
m=109   k non-p-rational for p=359  W=1 Hinfty=1
m=109   k non-p-rational for p=809  W=1 Hinfty=1
m=4061  k non-p-rational for p=2    W=1 Hinfty=1
m=4061  k non-p-rational for p=3    W=1 Hinfty=1
m=4061  k non-p-rational for p=5    W=1 Hinfty=1
m=4061  k non-p-rational for p=7    W=1 Hinfty=1
m=4061  k non-p-rational for p=11   W=1 Hinfty=1
m=4061  k non-p-rational for p=359  W=1 Hinfty=1
m=4061  k non-p-rational for p=433  W=1 Hinfty=1
m=4061  k non-p-rational for p=461  W=1 Hinfty=1
m=4061  k non-p-rational for p=2459 W=1 Hinfty=1
\end{verbatim}\ns

We have $\CW_k^\bp = \CH_k^\infty = 1$ except for $m = 42$
and $p = 3$, where $\order \CW_k^\bp = 3$. For $\CH_k^\infty \ne 1$, 
we have the following excerpts:

\smallskip
\ft\begin{verbatim}
m    p  k.cyc
79   3  [3],  142  3  [3],   223  3  [3],   229 3  [3], 235  3 [6], 254 3 [3],...
401  5  [5],  439  5  [5],   499  5  [5],   727 5  [5], 817  5 [5], 982 5 [5],...
577  7  [7],  1009 7  [7],   1087 7  [7],  1294 7  [7], 1534 7 [14],...
1297 11 [11], 4009 11 [11],  4762 11 [22], 5182 11 [11],...
4759 13 [13], 8101 13 [13],  8441 13 [13], 8647 13 [13],...
7054 17 [17],...
\end{verbatim}\ns

\subsubsection{Program using a ray class group}\label{P4bis}
Similarly, we give the usual program with ray class groups:

\smallskip
\ft\begin{verbatim}
{bp=2;Bp=10^4;for(m=2,10^3,if(core(m)!=m,next);rk=1;P=x^2-m;
k=bnfinit(P,1);forprime(p=bp,Bp,nup=2;if(p==2,nup=3);
kpnup=bnrinit(k,p^nup);Hpnup=kpnup.cyc;e=component(matsize(Hpnup),2);
R=0;for(j=1,e,c=Hpnup[e-j+1];if(Mod(c,p)==0,R=R+1));
if(R>rk,print("m=",m," Hk=",k.cyc," k is not ",p,"-rational"," rk(Tk)=",R-rk))))}
\end{verbatim}\ns

It finds again the previous examples, but takes a very long time:

${\sf time = 2min, 29,181\ ms\ up\ to\ 10^3}$.

\section{Program testing \texorpdfstring{$\order \CR_k(p) \ne 1$}{Lg} 
for any number field and non-bad \texorpdfstring{$p$}{Lg}}

\subsection{Matrix associated to \texorpdfstring{$\CR_k$}{Lg}}
Assume that $p$ is not ramified and does not divide the denominator
of the integer basis elements. Thus, $\{1, x, \ldots, x^{d-1}\}$ is
a $\Z_p$-basis.

\smallskip
Let $E_k :=  \{\varepsilon_1, \ldots ,\varepsilon_r \}$, $r = r_1+r_2 - 1$, be a 
fundamental system of units of $k$, congruent to $1$ modulo $(p)$. Put
$\frac{1}{p} \log(\varepsilon_\ell) = \sm_{j=0}^{d-1} b_j^\ell \cdot x^j$,
$\,\ell \in [1, r],\ \, b_j^\ell \in \Z_p$. 

\smallskip
Then we define the $(d, r)$-matrix:
\begin{equation}
\CM_k := 
 \begin{pmatrix}
b_0^1 & \!\!\!\! \ldots \!\!\!\! & b_0^r \\
 \vdots & \vdots &  \vdots  \\
b_{d-1}^1 & \!\!\!\! \ldots \!\!\!\! & b_{d-1}^r
\end{pmatrix}
\end{equation}

\noindent
where all the coefficients are considered modulo $p$ by the program. 
Then, the regulator $\CR_k(p) = \tor_{\Z_p} \big (\log(\CU_k)/
\log(\ov E_k) \big)$ is trivial if and only if the $\F_p$-rank of the 
matrix is maximal, equal to $r$;\,\footnote{\,Recall that if $\CR_k(p) = 1$ 
(as $\Z_p$-module), we say that $\CR_k(p)$ (as regulator) is invertible, 
which is more suitable since, to test the condition of invertibility, the 
program computes the $\F_p$-rank of the above matrix.}
the following part of the program 
was used in Program \ref{P1}, where $\log(\varepsilon_\ell)$ was 
computed in ${\sf Lell = LogE[ell]}$ (lift form of the logarithms).

\subsection{Two programs testing \texorpdfstring{$\CR_k(p) \ne 1$}{Lg}}
The following program can be used independently when $p \geq 2$ is not 
ramified and does not divide the denominators of the integer basis elements
(the other primes of Definition \ref{bad} are not bad for $\CR_k(p)$). 
Note that computation of units requests ${\sf bnfinit(P,1)}$).

\subsubsection{Program testing \texorpdfstring{$\CR_k(p) \ne 1$}{Lg} for 
non-bad primes}\label{P5}
Since the regulator is most often invertible, we only give an excerpt of
non-invertible examples:

\smallskip
\ft\begin{verbatim}
{P=x^4-17540*x^3-1542*x^2+4*x+1;d=poldegree(P);maxp=10^7;
print("P=",P);kP=nfinit(P);r1=kP[2][1];r2=kP[2][2];rE=r1+r2-1;
print("Z-rank of the unit group=",rE);k=bnfinit(P,1);Dk=k.disc;Zk=k.zk;
forprime(p=2,maxp,if(Mod(Dk,p)==0,print("bad prime ",p);next);
denB=1;for(j=1,d,Bzk=lift(Zk[j]);for(i=0,d-1,
denB=lcm(denominator(polcoeff(Bzk,i)),denB)));
if(Mod(denB,p)==0,print("bad prime ",p);next);
F=idealfactor(k,p);g=matsize(F)[1];M=1;
for(i=1,g,fi=component(F,1)[i][4];M=lcm(M,p^fi-1));
\\TEST ON THE REGULATOR:
N=2;if(p==2,N=3);mod=p^N;LogE=List;E=k.fu;
for(ell=1,rE,epsilonM=Mod(Mod(lift(E[ell]),mod),P)^M;
pi=Mod(lift(lift(epsilonM)),P)-1;if(pi==0,pi=Mod(mod,P));
\\COMPUTATION OF hatt:
hatt=1;
for(j=1,g,Primej=component(F,1)[j];
Deltaj=idealval(k,pi,Primej);rhoj= Primej[3];
hatx=floor((log(rhoj)-log(Deltaj*log(p)))/log(p))+1;
n=max(1,hatx);while(p^n*Deltaj-rhoj*(N+n)<0,n=n+1);
hatn=n;if(hatn==1,hattj=N);if(hatn>1,hattj=p^hatn;
while(hattj*Deltaj-(valuation(hattj,p)+N)*rhoj>0,hattj=hattj-1));
hatt=max(hatt,hattj));
Log=pi;t=1;pit=pi;while(t<=hatt,t=t+1;pit=pi*pit;
Log=Log+(-1)^(t+1)*pit/t);
\\REDUCTION OF Log MODULO mod:
Log=Mod(Mod(lift(Log),mod),P);listput(LogE,lift(lift(Log))));
\\THE MATRIX MLogE:
ME=List;for(ell=1,rE,Lell=LogE[ell];Lc=List;for(i=0,d-1,
coeff=Mod(polcoeff(Lell/p,i),p);listput(Lc,coeff));listput(ME,Lc));
MLogE=matrix(rE,d,u,v,ME[u][v]);RLog=matrank(MLogE);
\\if(RLog==rE,print("The ",p,"-regulator is invertible"));
if(RLog<rE,print("p=",p," rank(MLogE)=",RLog);
print("k is not ",p,"-rational: R non-invertible");
print("**********************")))}

x^4-17540*x^3-1542*x^2+4*x+1
bad prime 2
bad prime 3
p=5 rank(MLogE)=2
k is not 5-rational: R non-invertible
**********************
bad prime 11
bad prime 19
p=23 rank(MLogE)=1
k is not 23-rational: R non-invertible
**********************
p=29 rank(MLogE)=2
k is not 29-rational: R non-invertible
**********************
bad prime 43
p=47 rank(MLogE)=2
k is not 47-rational: R non-invertible
**********************
p=1873 rank(MLogE)=2
k is not 1873-rational: R non-invertible
**********************
p=2692343 rank(MLogE)=2
k is not 2692343-rational: R non-invertible
**********************
time = 15,410 ms for p<10^6.
time = 3min, 41,622 ms for p<10^7.
\end{verbatim}\ns

Up to $10^7$, all the other regulators are invertible.

\subsubsection{Program testing \texorpdfstring{$\CR_k(p) \ne 1$ for huge  
primes $p$, possibly bad}{Lg}}\label{P5bis}

It remains the case of a huge prime dividing the denominators of the 
coefficients of the integer basis elements. We assume that we have checked
that $\CH_k^\infty = \CW_k^\bp = 1$. The following example of an imaginary 
abelian field shows the method; this program finds again the previous examples,
but the method increase substantially the speed and must be used for primes 
such that ${\sf epsilonM=Mod(Mod(lift(E[ell]),mod),P)^M}$ is unworkable 
because of a huge ${\sf M}$ or $p$ being some bad prime. We consider, 
below, the polynomial ${\sf P=polsubcyclo(17*13*19,16)[1]}$.
The whole set of bad primes is:

\smallskip
\ft\begin{verbatim}
{2,3,5,7,11,13,17,19,151,331,727,773,1097,1109,1861,17489,19219,22441,33091,
 514081,1042381,1323871,28665829,127637651,3257494459,28161933524443}
\end{verbatim}\ns

\smallskip
The integer basis elements involve coefficients that are spectacularly large 
(numerators and denominators), and $p = 28161933524443$ divides the 
denominators.

\smallskip
Thus, the program is valid under the above assumptions; for small bad primes,
it can be cancelled by absurd calculations modulo $p$ due to denominators; in 
this example, only $13, 17,19$ (divisors of the discriminant) are bad primes. 
For all the others (except $p = 2$ for which $k$ is not $2$-rational), the program 
concludes that the regulator is invertible. For the list of $23$ primes we get only: 
${\sf time = 47,042\, ms}$.

\smallskip
Since the direct computation of $\pi = \varepsilon_i^M - 1 = 
\varepsilon_i^{28161933524443-1} -1$ is impossible (enormous time 
and memory), we use a classical induction with successive powers of 
${\sf m = 2}$, with reductions modulo the minimal polynomial ${\sf Z(x)}$ 
of the unit ${\sf epsilon}$ (of degree $\leq 8$ in this example), and 
reductions modulo ${\sf mod}$; having ${\sf x^M \equiv XM(x) 
\pmod {Z(x)}}$, we deduce ${\sf epsilonM = XM(epsilon)}$ since
${\sf Z(epsilon)=0}$; but the polynomial ${\sf XM(x)}$ is with small
degree and coefficients. The induction with powers of ${\sf m = 2}$ 
may be replaced by ${\sf m = 3},\,\ldots$ to have an optimal run:

\smallskip
\ft\begin{verbatim}
allocatemem(1000000000)
{P=polsubcyclo(17*13*19,16)[1];print("P=",P);d=poldegree(P);
\\p=28161933524443;
p=20988936657440586486151264256610222593863921;
kP=nfinit(P);r1=kP[2][1];r2=kP[2][2];
rE=r1+r2-1;k=bnfinit(P,1);F=idealfactor(k,p);g=matsize(F)[1];M=1;
for(i=1,g,fi=component(F,1)[i][4];M=lcm(M,p^fi-1));
N=2;mod=p^N;E=k.fu;LogE=List;
\\COMPUTATION OF pi for epsiloni:
for(i=1,rE,epsiloni=E[i];Z=minpoly(epsiloni);
m=2;XM=1;q=M;while(q>0,xM=x;
n0=floor(log(q)/log(m));
for(j=1,n0,xM=Mod(Mod(xM^m,mod),Z));
XM=lift(lift(xM*XM));q=q-m^n0);
epsiloniM=0;for(j=0,poldegree(XM),epsiloniM=
epsiloniM+polcoeff(XM,j)*epsiloni^j);
pi=Mod(Mod(lift(epsiloniM),mod),P)-1;
pi=Mod(lift(lift(pi)),P);if(pi==0,pi=Mod(mod,P));
\\COMPUTATION OF  hatt AND Log(epsiloni):
hatt=1;
for(j=1,g,Primej=component(F,1)[j];
Deltaj=idealval(k,pi,Primej);rhoj=Primej[3];
hatx=floor((log(rhoj)-log(Deltaj*log(p)))/log(p))+1;
n=max(1,hatx);while(p^n*Deltaj-rhoj*(N+n)<0,n=n+1);
hatn=n;if(hatn==1,hattj=N);if(hatn>1,hattj=p^hatn;
while(hattj*Deltaj-(valuation(hattj,p)+N)*rhoj>0,hattj=hattj-1));
hatt=max(hatt,hattj));Log=pi;t=1;pit=pi;while(t<=hatt,
t=t+1;pit=pi*pit;Log=Log+(-1)^(t+1)*pit/t);listput(LogE,lift(Log)));
\\MATRIX MLogE:
ME=List;for(ell=1,rE,Lell=LogE[ell]/p;Lc=List;
LellZ=Mod(nfalgtobasis(k,Lell),p);listput(ME,LellZ));
MLogE=matrix(rE,d,u,v,ME[u][v]);RLog=matrank(MLogE);
if(RLog==rE,print("The ",p,"-regulator is invertible"));
if(RLog<rE,print("p=",p," rank(MLogE)=",RLog);
print("k is not ",p,"-rational: R non-invertible"))}

P=x^16+x^15-352*x^14-259*x^13+27980*x^12-77359*x^11+2274382*x^10
    -17610397*x^9+219273240*x^8-1531270581*x^7+13716090512*x^6
    -27702275513*x^5+238906342429*x^4+658286108804*x^3
    -4878604426700*x^2-15240684713648*x+164805474108224
The 28161933524443-regulator is invertible
time = 2,112 ms.
\end{verbatim}\ns

For the prime number 
$p=20988936657440586486151264256610222593863921$ the
program finds the invertibility of ${\sf R}$ in ${\sf time = 36,255\ ms}$.

\section{Study of the kernel \texorpdfstring{$\order \CW_k^\bp(p) \ne 1$}{Lg} 
for any number field}
We will give three variants of a program about $\CW_k^\bp(p)$ (defined by the 
exact sequence \eqref{regulator} under Leopoldt's conjecture), for any $k$ 
and $p \leq [k : \Q]+1$. Since the condition $\CW_k^\bp(p) = 1$ is necessary 
for the $p$-rationality, it is useful to test this condition before any other one. 
Contrary to what one might think, except obvious cases as $p = 2$, $k$
quadratic, the determination of $\bigoplus_v \mu_p(k_v)$ is not immediate. 

\subsection{Local and global roots of unity with {\sc pari/gp}}
The instruction ${\sf nfrootsof1(k)[1]}$ gives the order of the global group
$\mu(k)$, but under the instruction ${\sf k=nfinit(P)}$ which may use many 
time since it computes unnecessary parameters for our purpose; to our 
knowledge, the case of the local groups $\mu_p(k_v)$ is not given and 
if we want to compute the decomposition of $p$ in $k$ (defined by 
${\sf P}$) or a cyclotomic extension $K$ of $k$ (defined by 
${\sf Q}$), this needs using ${\sf k=bnfinit(P)}$ and ${\sf K=bnfinit(Q)}$. 
We will also propose other methods.

\smallskip 
We assume that the field $k$ is given with the knowledge of a defining monic
irreducible polynomial ${\sf P}$ of degree ${\sf d}$. We propose three
programs depending on the order of magnitude of~${\sf d}$. For instance, in 
the last Program \ref{PW3}, for some polynomials, the sole instruction ${\sf nfinit(P)}$ 
would need: time = 29,075 ms and ${\sf nfinit(Q)}$ would be out of reach since
${\sf Q}$ is the compositum of ${\sf P}$ with various cyclotomic polynomials; 
so, we compute the local degrees from the $p$-adic factorisation of ${\sf P}$ 
and ${\sf Q}$, via ${\sf factorpadic(P,p,1)}$ and ${\sf factorpadic(Q,p,1)}$.

\smallskip
Nevertheless, in case of a problem, one may use another program; for instance,
with the polynomial ${\sf P=x^9-8*x^3-18*x+48}$ and ${\sf p=7}$, the first program
has a problem because of ${\sf k=nfinit(P)}$, probably with the unnecessary
computation of:
$${\sf nfdisc(P)=2165832547340333875200}=2^{30}\cdot 3^{11}
\cdot 5^2 \cdot 455461, $$ 

\noindent
while the two other programs need ${\sf time = 19 ms}$ and ${\sf time = 23 ms}$, respectively:

\smallskip
We know that $\CW_k^\bp(p) = 1$ for $p > d+1$; so the number of 
possible primes $p$ is finite and we will consider, in the programs, 
the set ${\sf PW}$ of primes $p \geq 2$ such that $p \leq d+1$. 

\smallskip
An obvious and frequent case giving $\CW_k^\bp(p) \ne 1$ is when 
$\mu_p(k) \ne 1$ (e.g., $p=2$) and when $p$ splits at least into two 
$p$-places in $k$.

\smallskip
So we get the following criterion depending on the number of $p$-places:

\smallskip
(i) $\mu_p(k) = 1$. Then, $\CW_k^\bp(p) \ne 1$ if and only if
there exists at least a $p$-place $v$ of $k$ such that $\mu_p(k_v) \ne 1$.

\smallskip
(ii) $\mu_p(k) \ne 1$. Then $\CW_k^\bp(p) \ne 1$ if and only if,
either $p$ splits at least into two $p$-places of $k$, or else $v \mid p$ is 
unique and $\order \mu_p(k_v) >  \order\mu_p(k)$.

\smallskip
To compute the local groups $\mu_p(k_v)$, we use the following criterion:

\begin{lemma}
Let $k$ be any number field defined via a polynomial $P$ of degree $d$.
Let $d_v := [k_v : \Q_p]$ for $v \mid p$. Let $p^n$ ($n \geq 0$ for $p > 2$,
$n \geq 1$ for $p = 2$) be the order of $\mu_p(k)$ and put $K=k(\mu_{p^{n+1}})$. 
Then locally, $\order \mu_p(k_v) > \order \mu_p(k)$ if and only if $v$ totally splits 
in $K/k$, whence, if and only if the local degree $d_w$ of any $w \mid v$ in $K/\Q$ 
is equal to $d_v$.
\end{lemma}

\begin{proof}
Obvious; note that $K/k$ is $p$-cyclic of non-trivial degree by definition of $p^n$.
\end{proof}

\subsection{Programs testing \texorpdfstring{$\CW_k^\bp \ne 1$
for all $p$}{Lg}}\label{PW}

In a computational point of view, we will first apply the test of decomposition
of $p$ in $k/\Q$ when $\mu_p(k) \ne 1$, then the tests of divisibilities of the 
$d_v$ by $p-1$, for at least one place. Whence the two possibilities:

\smallskip
(i) $\mu_p(k) = 1$. We test if there exists, at least, a non-trivial $\mu_p(k_v)$;

\smallskip
(ii) $\mu_p(k) \ne 1$ and $v \mid p$ unique (otherwise $\CW_k^\bp \ne 1$). 
We test if $\mu_p(k_v)$ is larger than $\mu_p(k)$.

\smallskip
To give more information, in case (ii) we test $\mu_p(k_v)$ even if 
$v$ is not unique (otherwise, put break instead of next); so it is
possible that $\order \mu_p(k_v) > \order \mu_p(k)$  will be computed
for several places.

\subsubsection{Program for small degrees}\label{PW1}

This program only uses ${\sf k=nfinit(P)}$ and ${\sf K=nfinit(Q)}$ for the
compositum ${\sf Q}$ of ${\sf P}$ with cyclotomic polynomials:

\smallskip
\ft\begin{verbatim}
PROGRAM FOR SMALL DEGREES:
{P=x^9-8*x^3-18*x+46;
\\P=x^9-8*x^3-18*x+48;
\\SET PW OF PRIMES p<=d+1:
d=poldegree(P);NW=0;PW=List;forprime(p=2,d+1,
listput(PW,p);NW=NW+1);print();print(P);print("PW=",PW);
\\COMPUTATION OF Muglobal:
k=nfinit(P);Muglobal=nfrootsof1(k)[1];
\\PRIME p IN PW AND LIST Ldv OF LOCAL DEGREES dv FOR k:
for(ell=1,NW,p=PW[ell];
F=idealfactor(k,p);g=matsize(F)[1];Ldv=List;
for(i=1,g,Primei=component(F,1)[i];
listput(Ldv,Primei[3]*Primei[4]));
\\CRITERION OF SPLITTING OF p:
wp=valuation(Muglobal,p);if(wp>0 & g>1,
print("W(",p,") non-trivial:wp>0 and g>1"," Ldv=",Ldv);
print("k is not ",p,"-rational");
print("**********************");next);
\\NOW ANY CASE MAY EXIST EXCEPT Muglobal>1 AND g>1.
\\TEST p-1 DIVIDES dv FOR AT LEAST ONE PLACE:
test=0;for(i=1,g,if(Mod(Ldv[i],p-1)==0,test=1;break));
if(test==0,print("p=",p,
" no local/global p-roots of unity");next);
\\LIST Ldw OF LOCAL DEGREES dw FOR K:
Q=polcompositum(P,polcyclo(p^(wp+1)))[1];K=nfinit(Q);
FKp=idealfactor(K,p);gK=matsize(FKp)[1];Ldw=List;
for(i=1,gK,PrimeiK=component(FKp,1)[i];
listput(Ldw,PrimeiK[3]*PrimeiK[4]));
\\TEST OF SPLITTING IN K/k:
Delta=poldegree(Q)/d;if(Delta!=p^wp*(p-1),
print("p=",p," no linear disjunction"));
print("p=",p," Mupglobal=",p^wp," Delta=",Delta);
\\LISTS OF LOCAL DEGREES IN ASCENDING ORDER:
Ldv=List(vecsort(Ldv));Ldw=List(vecsort(Ldw));
print("Ldv=",Ldv);print("Ldw=",Ldw);
\\COMPARISON OF THE LISTS, COMPUTATION OF dw/v:
t=0;deltav=0;while(t<g,t=t+1;delta=0;lambda=0;
while(delta<Delta,deltav=deltav+1;
lambdav=Ldw[deltav]/Ldv[t];print("dw/v(",p,")=",lambdav);
if(lambdav==1,print("Wv(",p,") non-trivial");
print("k is not ",p,"-rational");
print("**********************"));delta=delta+lambdav)))}

x^9-8*x^3-18*x+46
PW=List([2,3,5,7])
p=2 Mupglobal=2 Delta=2
Ldv=List([9])
Ldw=List([18])
dw/v(2)=2
p=3 Mupglobal=1 Delta=2
Ldv=List([3,6])
Ldw=List([6,12])
dw/v(3)=2
dw/v(3)=2
p=5 no local/global p-roots of unity
p=7 no local/global p-roots of unity
time = 88 ms.
\end{verbatim}\ns

\subsubsection{Program for medium degrees}\label{PW2}

Since the polynomial $Q$ varies and is of hight degree,
this program does not use ${\sf K=nfinit(Q)}$, replaced by the 
$p$-adic decomposition of $Q$; it only uses ${\sf k=nfinit(P)}$.

\smallskip
\ft\begin{verbatim}
PROGRAM FOR MEDIUM DEGREES:
{P=x^19+8*x^3-18*x-18;
\\P=x^9-8*x^3-18*x+48;
\\SET PW OF PRIMES p<=d+1:
d=poldegree(P);NW=0;PW=List;forprime(p=2,d+1,
listput(PW,p);NW=NW+1);print();print(P);print("PW=",PW);
\\COMPUTATION OF Muglobal:
k=nfinit(P);Muglobal=nfrootsof1(k)[1];
\\PRIME p:
for(ell=1,NW,p=PW[ell];
\\LIST Ldv OF LOCAL DEGREES dv FOR k:
F=idealfactor(k,p);g=matsize(F)[1];Ldv=List;
for(i=1,g,Primei=component(F,1)[i];
listput(Ldv,Primei[3]*Primei[4]));
\\CRITERION OF SPLITTING OF p:
wp=valuation(Muglobal,p);if(wp>0 & g>1,
print("W(",p,") non-trivial:wp>0 and g>1"," Ldv=",Ldv);
print("k is not ",p,"-rational");
print("**********************");next);
\\NOW ANY CASE MAY EXIST EXCEPT Muglobal>1 AND g>1.
\\TEST p-1 DIVIDES dv FOR AT LEAST ONE PLACE:
test=0;for(i=1,g,if(Mod(Ldv[i],p-1)==0,test=1;break));
if(test==0,print("p=",p,
" no local/global p-roots of unity");next);
\\LIST Ldw OF LOCAL DEGREES dw FOR K:
Q=polcompositum(P,polcyclo(p^(wp+1)))[1];Ldw=List;
FQ=factorpadic(Q,p,1);FQ=component(lift(FQ),1);
gK=matsize(FQ)[1];for(j=1,gK,listput(Ldw,poldegree(FQ[j])));
Delta=poldegree(Q)/d;if(Delta!=p^wp*(p-1),
print("p=",p," no linear disjunction"));
print("p=",p," Mupglobal=",p^wp," Delta=",Delta);
\\LISTS OF LOCAL DEGREES IN ASCENDING ORDER:
Ldv=List(vecsort(Ldv));Ldw=List(vecsort(Ldw));
print("Ldv=",Ldv);print("Ldw=",Ldw);
\\COMPARISON OF THE LISTS = COMPUTATION OF dw/v:
t=0;deltav=0;while(t<g,t=t+1;delta=0;lambda=0;
while(delta<Delta,deltav=deltav+1;
lambdav=Ldw[deltav]/Ldv[t];print("dw/v(",p,")=",lambdav);
if(lambdav==1,print("Wv(",p,") non-trivial");
print("k is not ",p,"-rational");
print("**********************"));delta=delta+lambdav)))}

x^19+8*x^3-18*x-18
PW=List([2,3,5,7,11,13,17,19])
p=2 Mupglobal=2 Delta=2
Ldv=List([19])
Ldw=List([38])
dw/v(2)=2
p=3 Mupglobal=1 Delta=2
Ldv=List([1,1,2,2,2,3,4,4])
Ldw=List([2,2,4,4,4,6,8,8])
dw/v(3)=2
dw/v(3)=2
dw/v(3)=2
dw/v(3)=2
dw/v(3)=2
dw/v(3)=2
dw/v(3)=2
dw/v(3)=2
p=5 no local/global p-roots of unity
p=7 no local/global p-roots of unity
p=11 no local/global p-roots of unity
p=13 no local/global p-roots of unity
p=17 no local/global p-roots of unity
p=19 no local/global p-roots of unity
time = 65 ms.
\end{verbatim}\ns

\subsubsection{Program for large degrees}\label{PW3}

This program does not use ${\sf k=nfinit(P)}$ nor ${\sf K=nfinit(Q)}$,
replaced by the $p$-adic factorizations of $P$ and $Q$.

\smallskip
\ft\begin{verbatim}
PROGRAM FOR VERY LARGE DEGREES:
allocatemem(800000000)
{P=x^29+8*x^3-18*x-18;
\\P=x^9-8*x^3-18*x+48;
\\SET PW OF PRIMES p<=d+1:
d=poldegree(P);NW=0;PW=List;forprime(p=2,d+1,
listput(PW,p);NW=NW+1);print();print(P);print("PW=",PW);
for(ell=1,NW,p=PW[ell];
\\COMPUTATION OF Muglobal(p):
w0=0;while(poldegree(polcompositum(P,polcyclo(p^w0))[1])==d,
w0=w0+1);wp=w0-1;
\\LIST Ldv OF LOCAL DEGREES dv FOR k:
Ldv=List;FP=factorpadic(P,p,1);FP=component(lift(FP),1);
g=matsize(FP)[1];for(j=1,g,dv=poldegree(FP[j]);listput(Ldv,dv));
\\CRITERION OF SPLITTING OF p:
if(wp>0 & g>1,print("W(",p,") non-trivial:wp>0 and g>1"," Ldv=",Ldv);
print("k is not ",p,"-rational");print("**********************");next);
\\NOW ANY CASE MAY EXIST EXCEPT Muglobal>1 AND g>1.
\\TEST p-1 DIVIDES dv FOR AT LEAST ONE PLACE:
test=0;for(i=1,g,if(Mod(Ldv[i],p-1)==0,test=1;break));
if(test==0,print("p=",p,
" no local/global p-roots of unity");next);
\\LIST Ldw OF LOCAL DEGREES dw FOR K:
Q=polcompositum(P,polcyclo(p^(wp+1)))[1];
Ldw=List;FQ=factorpadic(Q,p,1);FQ=component(lift(FQ),1);
gK=matsize(FQ)[1];for(j=1,gK,listput(Ldw,poldegree(FQ[j])));
\\TEST OF SPLITTING IN K/k:
Delta=poldegree(Q)/d;if(Delta!=p^wp*(p-1),
print("p=",p," no linear disjunction"));
print("p=",p," Mupglobal=",p^wp," Delta=",Delta);
\\LISTS OF LOCAL DEGREES IN ASCENDING ORDER:
Ldv=List(vecsort(Ldv));Ldw=List(vecsort(Ldw));
print("Ldv=",Ldv);print("Ldw=",Ldw);
\\COMPARISON OF THE LISTS = COMPUTATION OF dw/v:
t=0;deltav=0;while(t<g,t=t+1;delta=0;lambda=0;
while(delta<Delta,deltav=deltav+1;
lambdav=Ldw[deltav]/Ldv[t];
print("dw/v(",p,")=",lambdav);
if(lambdav==1,print("Wv(",p,") non-trivial");
print("k is not ",p,"-rational");
print("**********************"));delta=delta+lambdav)))}

x^29+8*x^3-18*x-18
PW=List([2,3,5,7,11,13,17,19,23,29])
p=2 Mupglobal=2 Delta=2
Ldv=List([29])
Ldw=List([58])
dw/v(2)=2
p=3 no local/global p-roots of unity
p=5 Mupglobal=1 Delta=4
Ldv=List([8,21])
Ldw=List([32,84])
dw/v(5)=4
dw/v(5)=4
p=7 no local/global p-roots of unity
p=11 no local/global p-roots of unity
p=13 no local/global p-roots of unity
p=17 Mupglobal=1 Delta=16
Ldv=List([3,4,6,16])
Ldw=List([48,64,96,256])
dw/v(17)=16
dw/v(17)=16
dw/v(17)=16
dw/v(17)=16
p=19 no local/global p-roots of unity
p=23 no local/global p-roots of unity
p=29 no local/global p-roots of unity
time = 10,683 ms.
\end{verbatim}\ns

Of course, the three programs give again immediately the obvious cases
of quadratic fields, a useful check. 

\smallskip
But for random polynomials, 
it is extremely difficult to find examples with $\CW_k^\bp(p) \ne 1$, 
except for $p = 2$ when for instance $2$ splits in $k$ since $-1 \in k$. 
The case $p = 3$ is often possible. We have rare examples with $p = 5$; 
the following ones show the complexity of the decomposition of the 
$p$-places and the structure of $\CW_k^\bp$; we use Program \ref{PW3}:

\smallskip
\ft\begin{verbatim}
{for(a=0,200,P=x^15+a*x^6-2*x+1;
if(polisirreducible(P)==0,next);
d=poldegree(P);NW=0;PW=List;forprime(p=2,d+1,
listput(PW,p);NW=NW+1);print();
print("a=",a);print(P);print("PW=",PW);
for(ell=1,NW,p=PW[ell];
w0=0;while(poldegree(polcompositum(P,polcyclo(p^w0))[1])==d,
w0=w0+1);wp=w0-1;
Ldv=List;FP=factorpadic(P,p,1);FP=component(lift(FP),1);
g=matsize(FP)[1];for(j=1,g,dv=poldegree(FP[j]);listput(Ldv,dv));
if(wp>0 & g>1,print("W(",p,") non-trivial:wp>0 and g>1"," Ldv=",Ldv);
print("k is not ",p,"-rational");print("**********************");next);
test=0;for(i=1,g,if(Mod(Ldv[i],p-1)==0,test=1;break));
if(test==0,print("p=",p,
" no local/global p-roots of unity");next);
Q=polcompositum(P,polcyclo(p^(wp+1)))[1];
Ldw=List;FQ=factorpadic(Q,p,1);FQ=component(lift(FQ),1);
gK=matsize(FQ)[1];for(j=1,gK,listput(Ldw,poldegree(FQ[j])));
Delta=poldegree(Q)/d;if(Delta!=p^wp*(p-1),
print("p=",p," no linear disjunction"));
print("p=",p," Mupglobal=",p^wp," Delta=",Delta);
Ldv=List(vecsort(Ldv));Ldw=List(vecsort(Ldw));
print("Ldv=",Ldv);print("Ldw=",Ldw);
t=0;deltav=0;while(t<g,t=t+1;delta=0;lambda=0;
while(delta<Delta,deltav=deltav+1;
lambdav=Ldw[deltav]/Ldv[t];
print("dw/v(",p,")=",lambdav);
if(lambdav==1,print("Wv(",p,") non-trivial");
print("k is not ",p,"-rational");
print("**********************"));delta=delta+lambdav))))}
time = 33,524 ms.

a=23
p=5 Mupglobal=1 Delta=4
Ldv=List([2,4,9])
Ldw=List([4,4,4,4,8,36])
dw/v(5)=2
dw/v(5)=2
dw/v(5)=1
Wv(5) non-trivial
k is not 5-rational
**********************
dw/v(5)=1
Wv(5) non-trivial
k is not 5-rational
**********************
dw/v(5)=2
dw/v(5)=4

a=48
p=5 Mupglobal=1 Delta=4
Ldv=List([2,4,9])
Ldw=List([4,4,4,4,4,4,36])
dw/v(5)=2
dw/v(5)=2
dw/v(5)=1
Wv(5) non-trivial
k is not 5-rational
**********************
dw/v(5)=1
Wv(5) non-trivial
k is not 5-rational
**********************
dw/v(5)=1
Wv(5) non-trivial
k is not 5-rational
**********************
dw/v(5)=1
Wv(5) non-trivial
k is not 5-rational
**********************
dw/v(5)=4

a=123
p=5 Mupglobal=1 Delta=4
Ldv=List([1,1,4,9])
Ldw=List([4,4,4,4,4,4,36])
dw/v(5)=4
dw/v(5)=4
dw/v(5)=1
Wv(5) non-trivial
k is not 5-rational
**********************
dw/v(5)=1
Wv(5) non-trivial
k is not 5-rational
**********************
dw/v(5)=1
Wv(5) non-trivial
k is not 5-rational
**********************
dw/v(5)=1
Wv(5) non-trivial
k is not 5-rational
**********************
dw/v(5)=4
\end{verbatim}\ns

For $a = 123$, the two places of residue degree $1$ in $k$ do not split
in $K_w/k_v$; the place of residue degree $4$ in $k$, totally splits in
$K_w/k_v$, so that $\mu_5$ is contained in $k_v$, then the last
place of $k$ does not split. To check, Program \ref{P2} finds:

\smallskip
\ft\begin{verbatim}
P=x^15+123*x^6-2*x+1  class group=[]
number of independent Z_p-extensions: 8
rk(Tk)=1 k is not 2-rational   Tk[2]=List([2])
rk(Tk)=1 k is not 5-rational   Tk[5]=List([5])
\end{verbatim}\ns

\section{Tests of \texorpdfstring{$p$}{Lg}-rationality in the literature}

We have decided against providing an extensive bibliography due to the vast 
number of articles on the subject; instead, we have limited ourselves to the 
seminal articles on the definition and characterization of the notion of 
$p$-rationality, which generally contain very comprehensive bibliographies, 
nevertheless we give papers describing some algorithms of computation.

\subsection{Overview of the main articles}
To our knowledge, the main papers giving some material to test the
$p$-rationality of a number field, are \cite{PV2015, Gra2019a, MR2019, 
HJ2026}. In \cite[Section 3]{PV2015} and \cite[Program II, p. 37]{Gra2019a},
the method is essentially based on the computation of ray class groups 
of modulus a sufficient power of $p$. As we have explained,
these computations use, for instance with {\sc pari/gp} package, the
instructions ${\sf k=bnfinit(P,1)}$ and ${\sf kpnup=bnrinit(k,p^{nup})}$, then
${\sf Hpnup=kpnup.cyc}$ giving the group structure, for ${\sf nup}$ large 
enough; but if $p$ is large, 
the programs are very slow and may take hours or days.
The authors of \cite{PV2015} do not provide an effective program; they 
only say ``All the computations have been done using the {\sc pari/gp} 
system''; they give interesting statistics for quadratic, cubic and quintic 
fields, with their interpretations via Cohen--Lenstra heuristics.

\smallskip
In \cite{MR2019} the examples about non $p$-rationality come from a deep
study of some isotopic components of the Galois group of the maximal 
$p$-ramified pro-$p$-extension of $k$, connected, in one way or another,
to the regulators since they use (as for many papers) cases where units 
are algebraically explicit.

\smallskip
The method of Hofmann--Johnston \cite{HJ2026} makes especially 
a focus on the $p$-adic regulator of the field, at the disadvantage 
of a general computation of the difficult degree $[\wt k H_k^\nr : \wt k] 
= \order \CH_k^\infty$; of course, when $k$ is totally real, 
$\wt k$ is the cyclotomic $\Z_p$-extension for which the 
computation of this index is obvious. The authors define the 
notion of quasi-$p$-rational field when $\CR_k = \CW_k^\bp = 1$.
Several applications are given about cyclotomic fields and constructions 
of Galois representations with open images in the spirit of Greenberg
methods.

\smallskip
Another fact is that many examples in the literature are such that the 
$p$-class group is trivial, so that $\CH_k^\infty = 1$.

\smallskip
Finally, in most cases, computing resources are either not accessible 
at all or only partially, and often require knowledge about the used package, 
without possibility of copy and past effective programs. Last but not least,
random polynomial, of degrees larger than $2$, $3$ or $4$, are not considered, 
even though doing so is necessary to understand $p$-rationality.

\subsection{Examples of abelian fields considered by Maire--Rougnant}
In \cite{MR2019} are given numerical results from the use of Balady-Lecacheux 
polynomials giving explicit units. Let's give (with Program \ref{P1}) few examples 
giving again results of \cite[Appendix A]{MR2019}:

\smallskip\noindent
{\bf (A1)} Cyclic cubic fields:
\ft\begin{verbatim}
P=x^3+309*x^2-10*x-1  class group=[3]
number of independent Z_p-extensions=1
Pbad=List([2,3,11,13,79])
p=3 rank of the 3-torsion group:Tk[p]=List([3])
k is not 3-rational
**********************
p=19 rank(MLogE)=1
k is not 19-rational
**********************
p=439 rank(MLogE)=1
k is not 439-rational
**********************
p=2887 rank(MLogE)=1
k is not 2887-rational
**********************
p=59060857 rank(MLogE)=1
k is not 59060857-rational
**********************
p=122648623 rank(MLogE)=1
k is not 122648623-rational
time = 10min, 46,093 ms.

P=x^3+12286733*x^2-9970*x-1  class group=[20811,3]
number of independent Z_p-extensions=1
Pbad=List([2,3,7,13,157,991,1051,9973])
p=3 rank of the 3-torsion group:Tk[p]=List([3,3,3])
k is not 3-rational
**********************
p=7 rank of the 7-torsion group:Tk[p]=List([7])
k is not 7-rational
**********************
p=991 rank of the 991-torsion group:Tk[p]=List([991])
k is not 991-rational
**********************
p=43 rank(MLogE)=1
k is not 43-rational
**********************
p=14783491 rank(MLogE)=1
k is not 14783491-rational
**********************
\end{verbatim}\ns

\smallskip\noindent
{\bf (A2)} Cyclic quartic field:
\ft\begin{verbatim}
x^4-7588*x^3-870*x^2+4*x+1  class group=[4]
number of independent Z_p-extensions=1
Pbad=List([2,3,5,11,73])
p=2 rank of the 2-torsion group:Tk[p]=List([4])
k is not 2-rational
**********************
p=17 rank(MLogE)=2
k is not 17-rational
**********************
p=19363829 rank(MLogE)=2
k is not 19363829-rational
**********************
p=26900513 rank(MLogE)=2
k is not 26900513-rational
**********************
\end{verbatim}\ns

\smallskip
These values are those obtained in \cite[Table A.2.1.]{MR2019}.
One may be surprised by the order of magnitude of these solutions
regarding our conjecture of finiteness of primes $p$ such that
$\CT_k(p) \ne 1$; but this is a phenomenon of $p$-adic property
of the units. Indeed, varying the field $k$ and fixing $p$ it is not 
difficult to get huge solutions; for instance set $m=1+10007^4$, 
then Program \ref{P4} gives immediately a non-$p$-rationality for
$p=10007$ and other large primes, but nothing else up to $10^9$; why ?

\smallskip
\ft\begin{verbatim}
m=10028029413722402 k non-p-rational for p=2 W=1 Hinfty=2^5
m=10028029413722402  k non-p-rational for p=10007
m=10028029413722402  k non-p-rational for p=503137
m=10028029413722402  k non-p-rational for p=21807883
\end{verbatim}\ns

\smallskip
One may force a finite number of $p$-adic properties of the units, but
this has no influence on the $q$-adic properties of the units for the other 
primes; moreover we do not know any way giving such ad-hoc solutions
{\bf with bounded discriminants}; a fundamental task would be to relate the 
$p$-adic regulator and the discriminant, which is only known for real 
framework (see \cite{Si1984} and subsequent papers on the subject); in
other words, we would need the $p$-adic version of this result, but this 
would be proving Leopoldt's conjecture.

\smallskip
We give some excerpts of the Balady--Washington polynomials considered in 
\cite[]{MR2019} (real cyclic fields of degree $4$), for $p \leq 10^7$:
$${\sf P=x^4+4*(s^3-s^2+2*s-1)*x^3-6*(s^2+1)*x^2+4*x+1, s \in [-20,20]}$$

\noindent
(Program \ref{P1} needing less that $10$\,mn for the above $41$
polynomials, $p$ up to $10^6$):

\smallskip\noindent
{\bf (B)} Cyclic quartic fields:
\ft\begin{verbatim}
P=x^4-33764*x^3-2406*x^2+4*x+1  class group=[4]
Pbad=List([2,3,5,67,107])
p=2 rank of the 2-torsion group Tk=1
k is not 2-rational
**********************
p=3 rank of the 3-torsion group Tk=1
k is not 3-rational
**********************
p=19 rank(MLogE)=2
k is not 19-rational
**********************
p=353 rank(MLogE)=2
k is not 353-rational
**********************
p=379 rank(MLogE)=2
k is not 379-rational
**********************
p=2273 rank(MLogE)=2
k is not 2273-rational
**********************
p=4159 rank(MLogE)=2
k is not 4159-rational
**********************
p=349051 rank(MLogE)=2
k is not 349051-rational
**********************
p=408809 rank(MLogE)=2
k is not 408809-rational
**********************

P=x^4+10300*x^3-1182*x^2+4*x+1  class group=[2,2,2]
Pbad=List([2,3,5,11,67])
p=2 rank of the 2-torsion group Tk=3
k is not 2-rational
**********************
p=3 rank of the 3-torsion group Tk=1
k is not 3-rational
**********************
p=29 rank(MLogE)=2
k is not 29-rational
**********************
p=43 rank(MLogE)=2
k is not 43-rational
**********************
p=73 rank(MLogE)=2
k is not 73-rational
**********************
p=271 rank(MLogE)=2
k is not 271-rational
**********************
p=409 rank(MLogE)=2
k is not 409-rational
**********************
p=1163 rank(MLogE)=2
k is not 1163-rational
**********************
p=9151 rank(MLogE)=2
k is not 9151-rational
**********************
p=28477 rank(MLogE)=2
k is not 28477-rational
**********************

P=x^4+12716*x^3-1356*x^2+4*x+1  class group=[2]
Pbad=List([2,3,5,227,619])
p=2 rank of the 2-torsion group Tk=1
k is not 2-rational
**********************

P=x^4-17540*x^3-1542*x^2+4*x+1  class group=[12,2,2]
Pbad=List([2,3,5,11,19,43])
p=2 rank of the 2-torsion group Tk=3
k is not 2-rational
**********************
p=3 rank of the 3-torsion group Tk=2
k is not 3-rational
**********************
p=5 rank of the 5-torsion group Tk=1
k is not 5-rational
**********************
p=23 rank(MLogE)=1
k is not 23-rational
**********************
p=29 rank(MLogE)=2
k is not 29-rational
**********************
p=47 rank(MLogE)=2
k is not 47-rational
**********************
p=1873 rank(MLogE)=2
k is not 1873-rational
**********************
p=2692343 rank(MLogE)=2
k is not 2692343-rational
**********************

P=x^4+18628*x^3-1740*x^2+4*x+1  class group=[4,4,2,2,2]
Pbad=List([2,3,5,11,73,97])
p=2 rank of the 2-torsion group Tk=5
k is not 2-rational
**********************
p=43 rank(MLogE)=2
k is not 43-rational
**********************
p=113 rank(MLogE)=2
k is not 113-rational
**********************
p=1021 rank(MLogE)=2
k is not 1021-rational
**********************
p=1667 rank(MLogE)=2
k is not 1667-rational
**********************
p=176677 rank(MLogE)=2
k is not 176677-rational
**********************
p=9697459 rank(MLogE)=2
k is not 9697459-rational
**********************
p=224035927 rank(MLogE)=2
k is not 224035927-rational
**********************
\end{verbatim}\ns

\subsection{Examples of abelian fields considered by Hofmann--Johnston}

Let's take, the first and last interesting numerical examples in \cite[Appendix A, 
Table 4]{HJ2026}, using polynomials giving totally real number fields $k$; we use 
Program \eqref{P1} for $p$ up to $10^7$:

\smallskip\noindent
{\bf (C)} Cyclic quintic field:
\ft\begin{verbatim}
P=x^5-68410*x^3-136820*x^2+1120453185*x-876769924  class group=[5]
number of independent Z_p-extensions=1
Pbad=List([2,3,5,7,179,6841])
**********************
p=2 Structure of the 2-torsion group:Tk=List([2,2,2,2])
k is not 2-rational
**********************
p=11 rank(MLogE)=3
k is not 11-rational
**********************
p=61 rank(MLogE)=3
k is not 61-rational
**********************
p=131 rank(MLogE)=3
k is not 131-rational
**********************
p=3701 rank(MLogE)=3
k is not 3701-rational
**********************
p=4421 rank(MLogE)=3
k is not 4421-rational
**********************
p=9311 rank(MLogE)=3
k is not 9311-rational
**********************
p=93581 rank(MLogE)=3
k is not 93581-rational
**********************
time = 4min, 55,199 ms.
\end{verbatim}\ns

\smallskip
Note that for $p = 5$, $\CH_k \ne 1$, but since the field $k$ is Galois real of 
degree $p$, the $\Log$-function is explicit (see \cite[Exercise IV.3.3.1]
{Gra2005}), and gives $\CT_k^{\Gal(k/\Q)} \simeq \CT_\Q = 1$.

\smallskip\noindent
{\bf (D)} Real abelian field of discriminant $401^{15}$
\ft\begin{verbatim}
P=x^20-3*x^19-119*x^18+177*x^17+5889*x^16-1381*x^15-152812*x^14-117981*x^13
+2164979*x^12+3491841*x^11-15663089*x^10-38883595*x^9+42546574*x^8
+187263905*x^7+52059719*x^6-314665312*x^5-318011368*x^4+4965810*x^3
+106222904*x^2+21505161*x-6072119  class group=[]
Pbad=List([2,3,5,7,11,13,17,19,83,401,487,797,1483,2897,11549,20903,
22277,24793,117329,771179,981263,2563217,5570717,81887591])
**********************
p=2 Structure of the 2-torsion group:Tk=List([2,2,2,2,2])
k is not 2-rational
**********************
p=3 Structure of the 3-torsion group:Tk=List([3])
k is not 3-rational
**********************
p=5 Structure of the 5-torsion group:Tk=List([5,5,5,5,5])
k is not 5-rational
**********************
p=17 Structure of the 17-torsion group:Tk=List([17])
k is not 17-rational
**********************
p=23 rank(MLogE)=18
k is not 23-rational
**********************
p=29 rank(MLogE)=18
k is not 29-rational
**********************
p=1093 rank(MLogE)=18
k is not 1093-rational
**********************
p=3251 rank(MLogE)=18p
k is not 3251-rational
**********************
p=19709 rank(MLogE)=17
k is not 19709-rational
**********************
p=38953 rank(MLogE)=18
k is not 38953-rational
**********************
p=228299 rank(MLogE)=17
k is not 228299-rational
**********************
p=272927 rank(MLogE)=18
k is not 272927-rational
**********************
p=1525021 rank(MLogE)=18
k is not 1525021-rational
**********************
p=1709287 rank(MLogE)=18
k is not 1709287-rational
**********************
time = 5h, 24min, 28,854 ms.
\end{verbatim}\ns

\smallskip
{\bf (E)} Another interesting polynomial, from \cite[Appendix A, Table 4]{HJ2026}, is:
$$P=x^8-9668*x^6+14526170*x^4-4580040976*x^2+224422007824, $$

\noindent
of quaternionic Galois group. Up to $10^7$, Program \ref{P1} gives:

\smallskip\noindent
\ft\begin{verbatim}
P=x^8-9668*x^6+14526170*x^4-4580040976*x^2+224422007824  class group=[2]
number of independent Z_p-extensions=1
Pbad=List([2,3,5,7,23,71,2417])
p=2 Structure of the 2-torsion group:Tk[p]=List([4,4])
k is not 2-rational
**********************
p=13 rank(MLogE)=6
k is not 13-rational
**********************
p=31 rank(MLogE)=6
k is not 31-rational
**********************
p=68059 rank(MLogE)=6
k is not 68059-rational
**********************
p=636863 rank(MLogE)=6
k is not 636863-rational
**********************
p=1546463 rank(MLogE)=6
k is not 1546463-rational
**********************
\end{verbatim}\ns

\section{Theoretical aspects and the \texorpdfstring{$\Log$}{Lg}-function}
The method is based, on the systematic use of the ordinary 
$p$-adic $\log$-function,
and in part, on the $\Log$-function, defined on $I_k \otimes \Z_p$
modulo $P_{k,\infty} := \bigcap_{n}\big( P_k(p^n) \otimes \Z_p\big)$, 
where $I_k$ is the group of prime-to-$p$ ideals of $k$ and $P_k(p^n)$
the ray group of modulus $(p^n)$; then, the kernel of $\Log$ is zero if and 
only if $k$ is $p$-rational (under Leopoldt's conjecture).

\subsection{Generalities on \texorpdfstring{$p$}{Lg}-ramification 
theory}\label{pramif}
Let $\wt k$ be the compositum of the $\Z_p$-extensions of $k$ and let
$H_k^\pr$ be the maximal abelian $p$-ramified pro-$p$ extension 
of $k$; then $\Gal(H_k^\pr/k) \simeq I_k \otimes \Z_p/P_{k,\infty}$.
Its torsion group $\CT_k = \Gal(H_k^\pr/\wt k)$ is finite, and a
characterization of the $p$-rationality is $\CT_k = 1$. Let $\CH_k$ be 
the $p$-class group of $k$ and let $H_k^\nr$ be the $p$-Hilbert 
class field of $k$. 

\smallskip
We will recall below that $\order \CT_k = \order \CH_k^\infty \times \order \CR_k 
\times \order \CW_k^\bp$, giving methods of computation.

\smallskip
The advantage of the $\Log$-function is that it gives the 
computation of the most difficult factor $\CH_k^\infty$, such that
$\CH_k/\CH_k^\infty$ is isomorphic to the Galois group of the 
maximal sub-extension of $H_k^\nr$ contained in $\wt k$.

\smallskip
Then, it is necessary to compute separately the two other invariants, $\CR_k$,  
the normalized $p$-adic regulator (\cite[Section 5]{Gra2018}) and $\CW_k^\bp 
= \bigoplus_v \mu_p(k_v) \big / \mu_p(k)$, in terms of local and global $p$-groups 
of roots of unity.

\smallskip
Most of the non-$p$-rational cases are due to small primes $p$ giving in 
particular $\CW_k^\bp(2) \ne 1$ and/or $\CW_k^\bp(3) \ne 1$; indeed, 
since $[\Q_p(\mu_p) : \Q_p] = [\Q(\mu_p) : \Q] = p-1$, a necessary 
condition to get $\CW_k^\bp(p) \ne 1$ is that there exists at least a $p$-place 
$v \mid p$ such that $\mu_p(k_v) \ne 1$, whence $p-1 \mid e_v$, the 
ramification index, giving $p \leq \max_{v \mid p}(e_v + 1)$ and a fortiori 
$p \leq d + 1$. 

\smallskip
Recall that $\CW_k^\bp(p)$ fixes the Bertrandias--Payan field $H_k^\bp$ 
and that $\CT_k^\bp := \Gal(H_k^\bp/\wt k)$ is probably the best invariant in 
$p$-ramification theory; cf. \cite[\S\,3.2]{Gra2014} defining weakly $p$-rational
fields when $\CT_k^\bp = 1$, under Leopoldt's conjecture.

\smallskip
The triviality of these three $\Z_p$-modules characterizes the $p$-rationality, 
according to the following diagram of abelian $p$-ramification theory over 
any number field $k$:
\unitlength=1.25cm 
\begin{equation*}
\begin{aligned}
\vbox{\hbox{\hspace{-4.0cm} 
\begin{picture}(10.0,3.3)
% horizontales
\put(6.4,2.50){\line(1,0){1.6}}
\put(8.8,2.50){\line(1,0){1.95}}
\put(3.7,2.50){\line(1,0){1.8}}
\put(4.3,2.24){\ft$\CH_k^\infty$}
\put(7.0,2.24){\ft$\CR_k$}
\put(9.35,2.24){\ft$\CW_k^\bp$}
\put(4.1,0.50){\line(1,0){1.35}} 
\put(1.6,0.50){\line(1,0){1.3}}
\put(1.3,0.40){$k$}
\bezier{300}(1.38,0.34)(3.5,0.0)(5.6,0.34)
\put(3.2,-0.1){\ft$\CH_k$}
\bezier{400}(3.7,2.7)(7.25,3.3)(10.8,2.7)
\bezier{400}(1.5,0.65)(7.25,1.0)(10.8,2.4)
\bezier{300}(5.9,2.3)(6.2,1.7)(10.85,2.4)
\put(6.1,1.75){\ft$\tor_{\Z_p}(\CU_k\!/\! \ov E_k)$}
\put(7.3,1.2){\ft$\CA_k$}
\put(7.0,3.14){\ft$\CT_k$}
% vertical
\put(3.5,0.7){\line(0,1){1.6}}
\put(5.7,0.7){\line(0,1){1.6}}
\bezier{350}(6.2,0.5)(8.5,0.6)(10.9,2.3)
\put(8.6,0.75){\ft$\CU_k\!/\! \ov E_k$}
\put(10.85,2.4){$H_k^\pr$}
\put(5.6,2.4){$\wt k H_k^\nr$}
\put(8.15,2.4){${H_k^\bp}$}
\put(3.35,2.4){$\wt k$}
\put(5.5,0.4){$H_k^\nr$}
\put(3.1,0.4){$\wt k \!\cap\! H_k^\nr$}
\end{picture}}} 
\end{aligned}
\end{equation*}
\unitlength=1.0cm

The group of $p$-principal local units $\CU_k$, whose images via local reciprocity 
maps give inertia groups in $H_k^\pr/H_k^\nr$, is of the form $\CU_k = 
\bigoplus_{v \mid p}\CU_k^v$, where $\CU_k^v := 1 + {\mathfrak p}_v$ 
(${\mathfrak p}_v =$ prime ideal associated to the place $v$). Then 
$\Gal(H_k^\pr/H_k^\nr) \simeq \CU_k/\ov E_k$, where $\ov E_k$ is the 
topological closure in $\CU_k$ of the diagonal embedding $\iota_p(E_k)$
of the group of $p$-principal global units of $k$. Thus, we have, under the 
Leopoldt conjecture:
\begin{equation}\label{regulator}
\left \{ \begin{aligned}
& \CR_k = \tor_{\Z_p} \big (\log(\CU_k)/\log(\ov E_k) \big), \\
& 1 \to \tor_{\Z_p}(\CU_k) \cdot \ov E_k/\ov E_k = \CW_k^\bp
\too \CU_k/\ov E_k \too \log(\CU_k)/\log(\ov E_k) \to 0,
\end{aligned}\right .
\end{equation}

\noindent
where $\log = \log_p = (\log_v)_{v \mid p}$ denotes the usual $p$-adic logarithm, 
with image in $\CO_k \big[\frac{1}{p} \big]$, where $\CO_k$ is the ring of integer
$\bigoplus_{v \mid p} \CO_k^v$ of $\bigoplus_{v \mid p} k_v$ with obvious notations. 
We will use the fact that if $p>2$ and if $v \mid p$ is not ramified, then
$\CU_k^v = 1 + p \cdot  \CO_k^v$,
$\log(\CU_k^v) = p \CO_k^v$ and $1 + p^2 \cdot  \CO_k^v = (\CU_k^v)^p$.

\smallskip
Then $k$ is $p$-rational if and only if $\CT_k = 1$. Thus, having checked the necessary 
conditions $\CW_k^\bp = \CH_k^\infty = 1$, formula \ref{regulator} reduces to the 
computation of $\CR_k = \tor_{\Z_p} \big(\log(\CU_k)/\log(\ov E_k) \big)$; in other
words, we do not need any ``regulator formula'' as that given by analytic formulas.
This is the object of Proposition \ref{matrix} in a computational point of view
using the fact that when $\CW_k^\bp = 1$, $\tor_{\Z_p}(\CU_k/\ov E_k) = 1$ 
is equivalent to $\log(\ov E_k)$ direct factor in $\log(\CU_k)$, which is easily 
characterized by the $\F_p$-rank of a matrix, especially since p is in principle 
unramified. In other words, as soon as we find $u \in \CU_k$ such that $u^p 
\in \ov E_k$, $u \notin \ov E_k$, $k$ is not $p$-rational.

\smallskip
Another important invariant is the Jaulent logarithmic class group $\wt \CH_k$
\cite{Jau1994, Jau1998, BJ2016} which measures the splitting in $\wt k/k^\cyc$, 
where $k^\cyc/k$ is the cyclotomic $\Z_p$-extension of $k$, so that an analogous 
schema does exist, replacing $H_k^\nr$ by the locally cyclotomic field $H_k^{\lc}$, 
then $\CH_k^\infty$ by $\wt \CH_k^\infty := \Gal(\wt k H_k^{\lc}/\wt k)$. 
The $p$-rationality and the logarithmic $p$-principality are not equivalent but 
$\CT_k = 1$ implies $H_k^{\lc} \subset \wt k$ since $H_k^{\lc} \subseteq H_k^{\pr}$
(all this under the Leopoldt and Gross-Kuz'min conjectures).

\begin{lemma}\label{rankH}
Under the Leopoldt conjecture for $p$ in the number field $k$ of degree $d$ and
signature $(r_1, r_2)$, with $d = r_1+ 2r_2$, a necessary condition
for the $p$-rationality of $k$ is that the $p$-rank of its $p$-class group $\CH_k$
be strictly less than $r_2+1$. If one knows that $k^\cyc/k$ is totally ramified, then 
the condition becomes: strictly less than $r_2$.
\end{lemma}

\begin{proof}
We have $\Gal(H_k^\pr/k) \simeq \Z_p^{r_2+1} \times \CT_k$; if $\CH_k^\infty = 1$, 
necessarily the $p$-Hilbert class field of $k$ is contained in $\wt k$, whence the claim.
\end{proof}

In cyclic $p$-extensions with sufficiently ramified primes, genus theory 
(Chevalley--Herbrand formula \cite[pp.\,402-406]{Che1933}) leads to large 
$p$-ranks, hence the non-$p$-rationality of $k$; at the opposite, this formula 
allows to characterize trivial $p$-class groups in cyclic $p$-extensions. 

\smallskip
See \cite{Gra1986} characterizing abelian $p$-rational $p$-extensions 
for $p = 2, 3$, via a fixed point formula for the $\Z_p$-modules $\CT$ 
in a $p$-extension, generalized in \cite[Theorem IV.3.3]{Gra2005}, which
introduces the notion of $p$-primitive ramification. The paper \cite{GJ1989} 
introduces the similar property of $p$-regularity. In \cite{Mov1988} are 
characterized the $p$-extensions in which $p$-rationality propagates; 
this uses the $p$-primitive ramification related to direct sums generated 
by suitable Frobenius elements $\big(\frac{\wt k/k}{{\mathfrak q}} \big)$ 
of images $\Log({\mathfrak q})$ in the isomorphism $\Gal(\wt k/k) \simeq 
\Log(I_k \otimes \Z_p)$ \cite[IV.3,~(b), Example IV.3.5.1]{Gra2005}.

\begin{proposition}\label{matrix}
Let $\varepsilon_1, \ldots , \varepsilon_r$ be a fundamental system of units of $k$.
Assume that $p \ne 2$ is unramified in $k/\Q$.
Put (with {\sc pari/gp} notations of the programs, where ${\sf Lell=LogE[ell]}$
denotes $\log(\varepsilon_\ell)$),
${\sf Lell/p = Mod(a_0^\ell + a_1^\ell*x + \cdots +
a_{d-1}^\ell* x^{d-1}, P)}, \ \,{\sf {a_j^\ell} \in \Z_p}$. 

\smallskip
Then, the normalized $p$-adic regulator $\CR_k$ is invertible if and only if the 
matrix ${\sf \big(Mod(a_j^\ell,p) \big)}$, $\ell \in [1, r]$, $j \in [0, d-1]$, is of maximal
$\F_p$-rank $r$.
\end{proposition}

\begin{proof} Under the assumptions, $\log(\CU_k) = p \cdot \CO_k$
and $\{1, x, \ldots, x^{d-1}\}$ plays the role of a local integer basis. 
We know that $\CR_k = \tor_{\Z_p} (\log(\CU_k)/\log(\ov E_k))$, where 
$\ov E_k$ is the topological closure in $\CU_k$ of the diagonal embedding 
$\iota_p(E_k)$ of the group of $p$-principal global units of $k$; then 
$\CR_k = 1$ if and only if $\log(\ov E_k)$ is a direct factor in $\log(\CU_k)$, 
whence the result.
\end{proof}

\subsection{The \texorpdfstring{$\Log$-function}{Lg}}\label{Logp}
The general study of $\CA_k := \Gal(H_k^\pr/k)$ and its torsion group
$\CT_k$ is done in our book, and we have the following canonical isomorphism, 
whatever the base-field $k$ and the prime $p$:
\begin{equation}\label{index}
\Gal(\wt k \cap H_k^\nr/k) \simeq \Log(I_k \otimes \Z_p)/
\Log(P_k \otimes \Z_p), 
\end{equation}

\noindent
where $I_k$ is the group of prime-to-$p$ ideals of $k$ and $P_k$ the sub-group 
of principal ideals; then: 
\begin{equation}\label{alpha}
\Log({\mathfrak a}) := \hbox{$\frac{1}{n}$} \log (\alpha) \!\!
\pmod {\Q_p \cdot \log (\ov E_k)}, \  \hbox{if ${\mathfrak a}^{n} = (\alpha)$, 
$(\alpha) \in P_k \otimes \Z_p$}. 
\end{equation}
 
Via the Artin Symbol, $\Log$ defines an isomorphism $\Log(I_k \otimes \Z_p) 
\simeq \Gal(\wt k/k)$. Whence the formula (product of three integers):
\begin{equation}\label{Tp}
\order \CT_k = \order \CH_k^\infty \cdot \order \CR_k \cdot \order\CW_k^\bp
= \frac{\order \CH_k}{\big[\Log(I_k \otimes \Z_p) : \Log(P_k \otimes \Z_p) \big]}
\cdot \order \CR_k \cdot \order \CW_k^\bp,
\end{equation}

\noindent
giving a test of non-$p$-rationality, the new needed invariant being the difficult 
denominator obtained from the $\Log$-function, which reduces to the usual
logarithm. 

\smallskip
The $\Q_p$-space $\Q_p \cdot \log (\ov E_k)$ is only relevant of 
the Leopoldt conjecture that we assume, and does not need numerical 
knowledge of the units (for instance, in the totally real case, $\Q_p\cdot 
\log (\ov E_k) \oplus \Q_p = \bigoplus_{v \mid p} k_v$, which explains
that this case does not really use $\Log$).

\begin{remarks}
(i) Set $\CH_k = \langle \Ccl({\mathfrak a}_1) \rangle \oplus \cdots 
\oplus \langle \Ccl({\mathfrak a}_h)\rangle$; the degree 
$[\wt k \cap H_k^\nr : k]$ only depends on the Fermat quotients of the 
$\alpha_i$'s deduced from the ${\mathfrak a}_i$'s, via the relations
\eqref{alpha}, ${\mathfrak a}_i^{n_i} = (\alpha_i)$.

\smallskip
(ii) In the particular case where $E_k \otimes \Z_p = 1$ (imaginary 
quadratic fields), the formula becomes $\order \CT_k = \ffrac{\order \CH_k}   
{\big[\log(I_k \otimes \Z_p) : \log(P_k \otimes \Z_p) \big]}
\cdot \order \CW_k^\bp$, with the usual $p$-adic logarithm.

\smallskip
(iii) If $k$ is a real quadratic field of fundamental unit $\varepsilon$, the formula
needs the computation of $[\wt k \cap H_k^\nr : k] = [k^\cyc \cap H_k^\nr : k]$; 
in that case, a direct computation is possible since the ramification in the 
cyclotomic $\Z_p$-extension $k^\cyc$ is well known.

\smallskip
But our purpose is to write a general program for fields of any signature; we 
begin with the case of imaginary quadratic fields (Program \ref{quadHinfty})
which reduces to the case of cyclic $p$-class groups for $p>2$ (Lemma 
\ref{rankH}). Then we give the general case (Program \ref{P6}). These two
programs use the $\Log$-function.
\end{remarks}

\section{Computational and programmatic considerations}

The difficulties, in a computational viewpoint, is the orders of magnitude of $p$ 
and/or the degree $d$ of the field $k$ defined via a monic irreducible polynomial 
$P \in \Z[x]$. So, if $p$ is small as well as $d$, our classical program using ray 
class groups of suitable modulus may be sufficient. However, under some 
circonstances, the {\sc pari/gp} computation of a ray class group with large 
modulus may be hopeless; in this case the new algorithm does not need ray
class groups and only computes $p$-adic logarithms of algebraic numbers. 

\smallskip
Another enormous difficulty is the order of magnitude of coefficients of some 
algebraic numbers (e.g., units ${\sf epsilon = Mod(e(x),P)}$, with possibly 
numerators of the coefficients of ${\sf e(x)}$ with thousands of digits), 
so that any banal operation takes much time. 
Since one only tests if $\order \CT_k \equiv 0 \pmod p$ or not, coefficients may be 
computed modulo a suitable modulus ${\sf mod = p^N}$ (${\sf N = 2}$ or $3$); 
then many reductions are done in the programs, which significantly increases the 
speed of execution.

\subsection{Bad primes}
In the $p$-adic $\log$ framework, a difficulty does appear when $p$ 
ramifies in~$k$ (slow convergence of the $\log$ maps) and when the 
algebraic numbers (units, pseudo-units, integer basis elements) have 
non-invertible denominators in their coefficients. So we may define the 
``bad primes'', few in number, as follows, and use for them, a ray class 
group of modulus $8$ for $p = 2$, $p^2$ for $p \geq 3$ not too large:

\begin{definition}\label{bad}
Let $k$ be a number field defined by the monic irreducible polynomial $P \in \Z[x]$, 
of degree $d$. Let $\CH_k$ be the $p$-class group of $k$ and let $\CD_k$ be its 
discriminant. Then we define, in the programs, the set ${\sf Pbad}$ of the following 
``bad'' prime numbers: 

\smallskip
(i) Primes $p \leq d+1$ ($p > d+1$ avoid to test $\CW_k^\bp(p) = 1$); set ${\sf PW}$.

\smallskip
(ii) Primes $p$ dividing $\CD_k$ ($p \nmid \CD_k$ 
simplifies logarithmic calculations (e.g., regulator $\CR_k$).

\smallskip
(iii) Primes $p$ dividing the class number $\order \BH_k$ ($\CH_k = 1$ 
avoid to test $\CH_k^\infty(p) = 1$).

\smallskip
(iv) Primes $p$ dividing the denominators of coefficients of the $d$ 
integer basis elements given under the form ${\sf k.zk[i] = Mod(ci(x),P)}$, 
where $c_i(x) \in \Q[x]$.
\end{definition}

All these ``bad primes'' are easily known from the very beginning. Then the 
programs process them first. However, in the case of quadratic fields 
$k = \Q(\sqrt{\pm m})$, we do not use ray class groups for these primes 
since for these fields, $\CW_k^\bp$ is known, the only possible 
denominator $2$ reduces bad primes to $p = 2$ (when $\pm m \equiv 1 
\pmod 4$), and some specificities can help; so it remains some more 
delicate computations of logarithms (see \S\,\ref{logalpha}): 

\smallskip
Ramified primes are given via ${\sf Dk=k.disc}$, then 
${\sf F=factor(Dk)}$, ${\sf component(F,1)}$ (list of the ${\sf matsize(F)[1]}$ 
prime divisors); the prime divisors of the class group are
computed in the same way from ${\sf k.no}$ and ${\sf k.cyc}$.

\smallskip
It will be necessary to compute separately $\CR_k$, which is the unique 
invariant that may be non-trivial for some $p \gg 0$; this is the goal of
Program \ref{P5bis}, after Program \ref{P5} for non-bad primes. We have 
conjectured, in various papers, that $\CR_k \ne 1$ happens only for a finite 
number of primes, in other words that $k$ is $p$-rational for almost all $p$.

\smallskip
Once we know that $\CW_k^\bp = \CR_k = 1$, we use the $\Log$-function
to compute the difficult invariant $\order (\CH_k/\CH_k^\infty) = 
\big[\Log(I_k \otimes \Z_p) : \Log(P_k \otimes \Z_p) \big]$ (see Section
\ref{Logp}, Formula \eqref{index}). In a computational aspect, the programs 
will use the {\sc pari/gp} instruction ${\sf bnfinit(P,1)}$ from which we get 
the class group ${\sf Hk=k.clgp}$ 
(then ${\sf Hk[1]}$ is the class number, ${\sf Hk[2]}$ gives the ${\sf h}$ 
components of the class group, and ${\sf Hk[3]}$ the representative 
ideals of the ${\sf h}$ components), and a system ${\sf Ek=k.fu}$ of 
${\sf rE = r1+r2-1}$ fundamental units.

\smallskip
Of course, the program using ray class groups is very short (few lines); moreover, 
with a suitable larger modulus, it computes the structure of $\CT_k$, but need
some interpretation of the result since {\sc pari/gp} gives the complete structure,
that is of the form $\Z_p^{r_2+1} \oplus \big[\oplus_{i=1}^s \Z/p^{t_i}\Z\big]$;
hence the $p$-rationality if and only if $s = 0$. So the reader may be surprised 
by the lengthy main Program \ref{P1}, but with an extremely fast 
running. This comes from a systematic optimization of some computations (e.g., 
large power of a large algebraic number and $\log$ computations with the usual 
series, that we justify in the next \S\,\ref{logalpha}).

\subsection{Principal algebraic numbers and prime-to-\texorpdfstring{$p$}{Lg} ideals}
If $\alpha$ is a prime-to-$p$ algebraic number of $k$, the computation of its $p$-adic
logarithm needs to work with a number congruent to $1$ modulo the product of the 
prime ideals of $k$ dividing $p$ (then called ``principal''); this obvious theoretical 
fact is not easy in a computational one since we must calculate a prime-to-$p$ 
power $M$ of $\alpha$ that may be extremely large.

\smallskip
For this, we must have all classes represented by prime-to-$p$ ideals 
${\mathfrak a}_i$; indeed, a suitable power ${\mathfrak a}_i^{n_i}$ gives
the pseudo-unit $\alpha_i$, needed for the computation
of $\Log({\mathfrak a}_i)$.
So one can use the following program giving the {\sc pari/gp} representative 
ideals ${\sf Ai}$ (from ${\sf Hk[3] = k.clgp[3]}$), but with ideals ${\sf Ai}$ non-necessarily 
prime to $p$. In the program, ${\sf Ni}$ is the norm of the ideal ${\sf Ai}$).
A list of prime-to-$p$ ideals ${\sf Bi}$ is given in ${\sf LHkp}$; the instruction 
${\sf z=idealappr(k,Aip)}$ is the key for this; recall its {\sc pari/gp} definition:

\smallskip
\ft\begin{verbatim}
idealappr(nf,x,{flag}): x being a fractional ideal, gives an element b such that 
v_p(b)=v_p(x) for all prime ideals p dividing x, and v_p(b)>=0 for all other p.

{P=x^2-2*3*5*7*11*17;p=23;k=bnfinit(P);Hk=k.clgp;LHkp=List([]);
dA=matsize(Hk[3])[2];for(j=1,dA,Ai=Hk[3][j];Ni=component(Ai,1)[1];
if(Mod(Ni,p)!=0,listput(LHkp,Ai));if(Mod(Ni,p)==0,Aip=idealmul(k,Ai,p);
z=idealappr(k,Aip);Bi=idealdiv(k,z,Aip);listput(LHkp,Bi)));
print("p=",p," Hk=",Hk," LHkp=",LHkp)}

p=23  Hk=[32,[4,2,2,2],[[29,27;0,1],[23,3;0,1],[85,0;0,1],[46,20;0,1]]]
LHkp=List([[29,27;0,1],[8535,15;0,5],[85,0;0,1],[27353261,5959649;0,1]])
\end{verbatim}\ns

The instruction ${\sf Aip=idealmul(k,Ai,p)}$ does not modify the class and
makes certain that ${\sf idealdiv(k,z,Aip)}$ will be prime to {\it all $p$-places}
when $p$ splits. Here, $23$ splits in $k$; the two ideals of norms $23$  and 
$46$ have been replaced by ideals of norms $42675$ and $27353261$, 
respectively.

\subsection{Relation \texorpdfstring{$\alpha^M \equiv 1 \pmod 
{\prod_{\mathfrak p \mid p} {\mathfrak p}}$}{Lg}}
We define a minimal integer $M$ (depending on the residue degrees of the 
$p$-places) such that for any prime-to-$p$ number $\alpha \in k^\times$, 
$\alpha^M$ is congruent to $1$ modulo all the $p$-places of $k$. In the 
programs, $\pi := \alpha^M - 1$ writes ${\sf pi = alphaM - 1}$. It generates 
an ideal of the form:
$$(\pi) = {\mathfrak a} \cdot \hbox{$\prod_{\mathfrak p \mid p}$}
\ {\mathfrak p}^{1 + \delta_i(\alpha)}, \ \hbox{$\delta_i(\alpha) \geq 0$, 
and ${\mathfrak a}$ prime to $p$}.,$$

\noindent
where the ${\mathfrak p}^{\delta_i(\alpha)}$'s play the role of ${\mathfrak p}$-Fermat 
quotients of $\alpha$, most often trivial.
If $p$ is unramified in $k/\Q$, then $(\pi) = {\mathfrak a} \cdot (p) \cdot
{\prod_{\mathfrak p \mid p} {\mathfrak p}^{\delta_i(\alpha)}}$.

\subsubsection{Computation of \texorpdfstring{$M$}{Lg}}\label{algo0}
We assume that the residue degrees of the $p$-places of $k$ are
given in the list ${\sf Lf=([f1, \ldots, fg])}$. A canonical solution
is the least common multiple $M$ of the ${\sf p^{fi}-1}$.
But $M$ can be so large that one must compute modulo
${\sf mod = p^N}$ some coefficients.

\subsubsection{Computation of 
\texorpdfstring{$\alpha^M \pmod{p^N}$}{Lg}}\label{algo1}

For instance, the following program, with a large unramified prime $p$,
computes the prime ideals ${\sf Primei}$ dividing $p$ (list ${\sf Ideal}$), 
ramification indices (list ${\sf Lrho}$), residue degrees (list ${\sf Lf}$), 
the integer ${\sf M}$, and ${\sf \alpha^M \pmod{p^N}}$:

\smallskip
\ft\begin{verbatim}
allocatemem(800000000)
{P=x^17-31*x-1;p=100000007;N=2;mod=p^N;alpha=Mod(x,P);
k=bnfinit(P);Ideal=List;Lrho=List;Lf=List;F=idealfactor(k,p);
g=matsize(F)[1];for(i=1,g,Primei=component(F,1)[i];listput(Ideal,Primei);
rhoi=Primei[3];fi=Primei[4];listput(Lrho,rhoi);listput(Lf,fi));
M=1;for(i=1,g,Mi=p^Lf[i]-1;M=lcm(M,Mi));
alphaM=Mod(Mod(lift(alpha),mod),P)^M;
print("Lrho=",Lrho," Lf=",Lf);print("M=",M);
print("alphaM=",lift(lift(alphaM)))}

Lrho=List([1,1,1,1,1]) Lf=List([1,1,2,4,9])
M=100000085000033120007822601247358141463991700
                 5539111378155216975637693591518786629056793341442400
                 
alphaM=3409200238644*x^16+5484595083921629*x^15+2119742948381996*x^14
+8987355229114822*x^13+7850508049535525*x^12+7732600941282028*x^11
+7477157623400997*x^10+2088556446198941*x^9+1916618434163281*x^8
+6251496237604706*x^7+400278228019474*x^6+326693722868559*x^5
+8452449591671430*x^4+6869203380844203*x^3+9922633494584296*x^2
+439241030746870*x+6959356087154893
=1+100000007*(34092*x^16+54845947*x^15+21197428*x^14+89873546*x^13
+78505075*x^12+77326004*x^11+74771571*x^10+20885563*x^9+19166183*x^8
+62514958*x^7+4002782*x^6+3266937*x^5+84524490*x^4+68692029*x^3
+99226328*x^2+4392410*x+69593556)
time = 1,901 ms.
\end{verbatim}\ns

\noindent
The run takes ${\sf time = 1,901\, ms}$, the main part being due to the instruction 
${\sf k=bnfinit(P)}$ for $P$ of large degree. It is clear that, for huge primes $p$, the 
discriminant of $k$ and its class number are in general prime to $p$, so that we will 
have the $p$-adic approximations modulo ${\sf mod = p^2}$ (e.g., for logarithms).

\subsection{Computation of \texorpdfstring{$\log(\alpha^M)$}{Lg} modulo
\texorpdfstring{$p^N$}{Lg}}
\label{logalpha}
Set $\alpha^M = 1+ \pi$, where $(\pi) =: {\mathfrak a} \cdot \prod_{\mathfrak p 
\mid p} {\mathfrak p}^{\Delta_{\mathfrak p}}$ with $\Delta_{\mathfrak p} = 
1 + \delta_{\mathfrak p}$. 
As usual, we must estimate a minimal number of terms of the series $\log(\alpha^M) 
= \sum_{t=1}^\infty (-1)^{t+1} \cdot \pi^t/t$ that may give a congruence, of the form:
$$\log(\alpha^M) = \log(1 + \pi) \equiv \sm_{t=1}^{\hat t} (-1)^{t+1} \cdot \pi^t/t 
\pmod{p^N\! \cdot \CO_k}, $$
 
\noindent
for $N$ equal in general to $2$ or $3$ depending on $p$.

\subsubsection{Minimal index $\hat t$ of the $\log$ series}
Let $g$ be the number of prime ideals ${\mathfrak p} \mid p$ of $k$; let 
$\rho_{\mathfrak p} \geq 1$ be the ramification index of ${\mathfrak p}$ 
in $k/\Q$. Then the relation 
$(p) = \prod_{{\mathfrak p} \mid p} {\mathfrak p}^{\rho_{\mathfrak p}}$ implies 
that a computation modulo $p^N$ is satisfied as soon as it is satisfied modulo 
each power ${\mathfrak p}^{N \cdot \rho_{\mathfrak p}}$.

\smallskip
Since the ${\mathfrak p}$-valuations $\Delta_{\mathfrak p}$ of the number $\pi$ 
can take any value $\geq 1$, the computation of the minimal limit term of the 
series strongly depends on the ramification indices $\rho_{\mathfrak p}$ of 
$p$ in $k$, then is not so obvious and must work in any case in the programs; 
in particular we can use majorations of $\hat t$ at will. We fix, once for all, 
$N = 2$ if $p >2$, $N = 3$ if $p = 2$:

\begin{proposition}\label{prop1}
Let $\alpha = 1 + \pi$, where $(\pi) =  {\mathfrak a} \cdot 
\prod_{\mathfrak p \mid p} {\mathfrak p}^{\Delta_{\mathfrak p}}$, 
${\mathfrak a}$ prime-to-$p$ ideal,
$\Delta_{\mathfrak p} \geq 1$ for all ${\mathfrak p} \mid p$. Put
$x_{\mathfrak p} := \big[\log(\rho_{\mathfrak p}) - \log(\Delta_{\mathfrak p} 
\cdot \log(p)) \big]/\log(p)$ (here, $\log = \log_\infty)$. 

\smallskip
Let $\hat n_{\mathfrak p}$ be the minimal 
integer $n$ such that: $n \geq \max(1, [x_{\mathfrak p}]+1)$
and $p^n \cdot \Delta_{\mathfrak p} - n \cdot \rho_{\mathfrak p} 
\geq N \cdot \rho_{\mathfrak p}$. 

\smallskip
Let $\hat t_{\mathfrak p}$ be the maximal integer 
$t \in[1, p^{\hat n_{\mathfrak p}} - 1]$ such that:
$p^t \cdot \Delta_{\mathfrak p} - v_p(t) \cdot \rho_{\mathfrak p} 
< N \cdot \rho_{\mathfrak p}$).

\smallskip
a) If $\hat n_{\mathfrak p} = 1$ (which supposes $p \geq (N+1) \cdot 
\rho_{\mathfrak p}/\Delta_{\mathfrak p}$), then $\hat t_{\mathfrak p} = N$
 is suitable.
 
b) Let $\hat t := \max_{{\mathfrak p} \mid p} (\hat t_{\mathfrak p}$).
Then $\sum_{t=1}^\infty (-1)^{t+1} \cdot \pi^t/t \equiv 
\sum_{t=1}^{\hat t} (-1)^{t+1} \cdot \pi^t/t \pmod{p^N}$.  

\smallskip
c) The number $\log(1+\pi)$ is in $\CO_k$ as soon as 
$\rho_{\mathfrak p} \leq 2$ for all ${\mathfrak p} \mid p$;
then $\log(1+\pi) \in \pi \cdot \CO_k$ as soon as 
$\rho_{\mathfrak p} \leq p-1$ for all ${\mathfrak p} \mid p$.
\end{proposition}

We note that if $\hat n_{\mathfrak p} = 1$, then $p \geq (N+1) 
\cdot \rho_{\mathfrak p}/\Delta_{\mathfrak p}$. In the frequent case 
$\Delta_{\mathfrak p} = \rho_{\mathfrak p} = 1$, this applies only 
for $p \geq 3$. For $p = 2$, this applies if $\rho_{\mathfrak p} = 1$ 
and if moreover $\Delta_{\mathfrak p} \geq 2$. Then 
$\hat t_{\mathfrak p}$ is the maximal integer $t \in [1, p-1]$ 
such that $p^t  < N \cdot \rho_{\mathfrak p}/\Delta_{\mathfrak p}$.

\begin{proof}
The purpose is to consider the following steps for each prime ideal ${\mathfrak p}$:

\smallskip
(i) Determination of $\hat n_{\mathfrak p} \geq 1$, then proof that, in the
interval $[p^{\hat n_{\mathfrak p}}, \infty]$, the terms of the form 
$\pi^{p^n}/p^n$, $n \geq \hat n_{\mathfrak p}$, have a strictly 
increasing ${\mathfrak p}$-valuation and can be neglected modulo $p^N$; 

\smallskip
(ii) Proof that $\pi^n/n \equiv 0 \pmod {p^N}$ for all $n \in 
[p^{\hat n_{\mathfrak p}}, p^{\hat n_{\mathfrak p}+1}]$, for any 
value of $\hat n_{\mathfrak p} \geq 1$.

\smallskip
(iii) Determination of $\hat t_{\mathfrak p}$ in $[p^{\hat n_{\mathfrak p}-1}, 
p^{\hat n_{\mathfrak p}}[$.

\smallskip
Let $g(x) :=  p^x \cdot \Delta_{\mathfrak p} - x \cdot \rho_{\mathfrak p}$, 
for $x \geq 0$ (${\mathfrak p}$-valuation function of the terms $\pi^{p^x}/p^x$); 
its derivative is $g'(x) = p^x \cdot \Delta_{\mathfrak p} 
\cdot \log(p) - \rho_{\mathfrak p}$, giving the minimum of $g(x)$ for:
$$x_{\mathfrak p} = \ffrac{\log(\rho_{\mathfrak p}) - \log(\Delta_{\mathfrak p} 
\cdot \log(p))}{\log(p)}. $$ 

Thus, $g(x)$ is strictly increasing from $x = x_{\mathfrak p}$. One verifies 
that when ${\mathfrak p}$ does not ramify, then $g(x_{\mathfrak p}) = 
(1+\log(\Delta_{\mathfrak p} \cdot \log(p))/\log(p) >0$ 
for all $p\geq 2$; this is not true if ${\mathfrak p}$ ramifies (e.g., 
$\Delta_{\mathfrak p} = 1$ and $\log(\rho_{\mathfrak p}) > 1+\log(\log(p))$). 

\smallskip
So, there exists $\hat n_{\mathfrak p} \geq \max(1, [x_{\mathfrak p}]+1)$ 
such that $v_{\mathfrak p}(\pi^{p^n}/p^n) \geq N \cdot \rho_{\mathfrak p}$ 
for all $n \geq n_{\mathfrak p}$, thus: 
\begin{equation}\label{pn}
v_{\mathfrak p}(\pi^{p^n}/p^n) =
p^n \cdot \Delta_{\mathfrak p} - n \cdot \rho_{\mathfrak p} 
- N \cdot \rho_{\mathfrak p} \geq 0, \  \hbox{for all $n \geq n_{\mathfrak p}$}. 
\end{equation}

If $x_{\mathfrak p} \leq 1$, then $g(x)$ is increasing and positive for all $x \geq 1$;
this condition is equivalent to $\rho_{\mathfrak p} \leq p \cdot \Delta_{\mathfrak p}
\cdot \log(p)$; it is satisfied for all unramifief $p \geq 2$; for quadratic fields
it is satisfied for $p \geq 3$ whatever the ramification.

\smallskip
So, it remains to prove that all $t \in ]p^n, p^{n+1}[$,
$n \geq \hat n_{\mathfrak p}$, may be neglected modulo $p^N$.
For this, it suffices to prove the inequality $v_{\mathfrak p}(\pi^{t}/t) \geq 
N \cdot \rho_{\mathfrak p}$ in this interval; put
$t = p^n + \lambda p^r$, for $\lambda \in [1, p-1]$
and $r \in [0, n]$ ($\lambda \ne p-1$ if $r = n$); we have, using \eqref{pn}:
\begin{equation}
\left \{ \begin{aligned}
v_{\mathfrak p}(\pi^{t}/t) & = t \cdot \Delta_{\mathfrak p} - r \cdot \rho_{\mathfrak p} =
(p^n + \lambda p^r) \cdot \Delta_{\mathfrak p} - r\cdot \rho_{\mathfrak p} \\
& = p^n \cdot \Delta_{\mathfrak p} + \lambda p^r \cdot \Delta_{\mathfrak p} 
- r \cdot \rho_{\mathfrak p}
\geq n\cdot \rho_{\mathfrak p} + N \cdot \rho_{\mathfrak p} +  
\lambda p^r \cdot \Delta_{\mathfrak p} - r\cdot \rho_{\mathfrak p} \\
& \geq N \cdot \rho_{\mathfrak p} + (n - r)\cdot \rho_{\mathfrak p}  
+ \lambda p^r \cdot \Delta_{\mathfrak p} > N \cdot \rho_{\mathfrak p}.
\end{aligned} \right.
\end{equation}

Whence the claims; the computations of $\hat n_{\mathfrak p}$
and $\hat t_{\mathfrak p}$ being solved by the program using the 
instruction ${\sf while}$; it will then take $\hat t = \max_{{\mathfrak p} \mid p}
(\hat t_{\mathfrak p})$ obtained from the above computations.

\smallskip
Assume that $\rho_{\mathfrak p} \leq 2$ for all ${\mathfrak p} \mid p$ 
and put $w := v_p(t)$ for $t \geq 1$; then: 
$$v_{\mathfrak p}(\pi^t/t) = t \cdot \Delta_{\mathfrak p} - w \cdot \rho_{\mathfrak p}
\geq p^w \Delta_{\mathfrak p}
- w \cdot \rho_{\mathfrak p} \geq p^w -  2 w, $$ 

\noindent
which is non-negative for all $w \geq 0$ and all $p \geq 2$; 
for the second claim, consider $\log((1+\pi)^p)$. Whence the
property (c), which applies in particular to quadratic fields.
\end{proof}

\begin{remark}\label{decrease}
For the speed of program execution, we are obliged to take all values of $t$
in a suitable interval $[1, \hat t]$, but many terms $\pi^t/t$ of large
valuation are useless. This explains that we must define $\hat t$ 
taking first the value $p^{\hat n}$, then decrease it as long as the 
valuations remain too large; the final value is not always $p^{\hat n - 1}$. 
For instance, for $p = 2$, $\Delta = 1$, $\rho = 7$, $N=3$:

\smallskip
\ft\begin{verbatim}
{Delta=1;rho=7;N=3;n=1;
while(2^n*Delta-(N+n)*rho<0,print("n=",n);n=n+1);hatn=n;print("hatn=",hatn);
t=2^hatn + 20;while(t*Delta-(valuation(t,2)+N)*rho>0,print("t=",t);t=t-1);print("hatt=",t);
for(t=1,100,print(t*Delta-(valuation(t,2)+N)*rho))}

n=1, n=2, n=3, n=4, n=5,  hatn=6 
t=64, t=63, t=62, t=61, t=60, t=59, t=58, t=57, 
t=56, t=55, t=54, t=53, t=52, t=51, t=50,t=49,
hatt=48     
Indeed:
values of t*Delta-(valuation(t,2)+N)*rho, from t=1 to t=48, then beyond t=49:
-20,-26,-18,-31,-16,-22,-14,-34,-12,-18,-10,-23,-8,-14,-6,-33,-4,-10,-2,-15,
0,-6,2,-18,4,-2,6,-7,8,2,10,-24,12,6,14,1,16,10,18,-2,20,14,22,9,24,18,26,-1,

28,22,30,17,32,26,34,14,36,30,38,25,40,34,42,1,44,...
\end{verbatim}\ns

Many terms are useless, but this does not affect the speed of the computations.
\end{remark}

Let's give the basic algorithm for a fixed prime ideal ${\mathfrak p}$, 
where ${\sf Delta}$ is the ${\mathfrak p}$-valuation of ${\sf pi}$ and 
${\sf rho}$ the ramification index of ${\mathfrak p}$:

\subsubsection{Computation of \texorpdfstring{$\hat t_{\mathfrak p}$
and of the $\Log$ series up to 
$\pi^{\hat t_{\mathfrak p}}/\hat t_{\mathfrak p}$}{Lg}}\label{algo2}

We give fictive data for which some coefficients of ${\sf Log}$ have
the denominator ${\sf p=2}$ and where the ramification index is ${\sf rho >2}$; 
so, the reduction modulo ${\sf mod = p^N}$ gives some difficulties with 
the instructions ${\sf Mod(Log,p)}$; the following program is valid for any 
data ${\sf p}$, ${\sf P}$, ${\sf Delta}$, ${\sf rho}$, ${\sf N}$, ${\sf mod=p^N}$, 
${\sf pi}$, obtained from the main programs testing the ${\sf p}$-rationality:

\smallskip 
\ft\begin{verbatim}
{p=2;P=x^8-2;Delta=1;rho=8;N=3;mod=p^N;pi=Mod(x,P);
\\COMPUTATION OF THE MINIMUM hatx OF g(x)=p^x.Delta-x.rho:
hatx=floor((log(rho)-log(Delta*log(p)))/log(p))+1;
\\COMPUTATION OF THE MAXIMAL INDICE hatn:
n=max(1,hatx);while(p^n*Delta-(N+n)*rho<0,n=n+1);
hatn=n;
\\COMPUTATION OF hatt FOR A PRIME:
if(hatn==1,hatt=N);if(hatn>1,
hatt=p^hatn;while(hatt*Delta-(valuation(hatt,p)+N)*rho>0,hatt=hatt-1));
\\COMPUTATION OF Log MOD mod:
Log=pi;t=1;pit=pi;while(t<=hatt,t=t+1;pit=pi*pit;Log=Log+(-1)^(t+1)*pit/t;
\\VALUATION OF THE DENOMINATOR OF Log:
vLog=min(valuation(lift(Log),p),0);Log=lift(Log)*p^-vLog;
\\REDUCTION MODULO mod*p^-vLog AND RESTITUTION OF Log:
Log=lift(Mod(Log,mod*p^-vLog));Log=Mod(Log*p^vLog,P));
print("n=",n," hatn=",hatn," hatt=",hatt);print("Log=",Log);
print("val(norm(Log))=",valuation(norm(Log),p))}

n=7 hatn=7 hatt=64
Log=Mod(x^7+3/2*x^6+3*x^5+13/4*x^4+5*x^3+9/2*x^2+7*x,x^8-2)
valuation(norm(Log),2)=-12
time = 1 ms.
\end{verbatim}\ns

Note that if $p$ is unramified such that $p \cdot \Delta_{\mathfrak p} - (1+N) > 0$, 
for all ${\mathfrak p} \mid p$, then $\hat x < 0$ and $\hat t = 2$ is suitable, 
without more computations; this is the case when $p \geq [k : \Q] \cdot (N+1)$.

\smallskip
In the opposite case, a tricky example may be when $p=2$ with large
ramification indices and the number $\pi$ such that $\Delta_{\mathfrak p} = 1$
for all ${\mathfrak p} \mid p$; we give an example in which the
numerical data may be modified at will, taking a unit:

\ft\begin{verbatim}
allocatemem(800000000)
{p=2;P=x^16-32*x^4+8*x-2;N=3;mod=p^N;k=bnfinit(P,1);alpha=k.fu[1];
\\DECOMPOSITION OF p AND COMPUTATION OF M:
Ideal=List;Lrho=List;Lf=List;M=1;F=idealfactor(k,p);g=matsize(F)[1];
for(i=1,g,Primei=component(F,1)[i];
listput(Ideal,Primei);rhoi=Primei[3];listput(Lrho,rhoi);
fi=Primei[4];M=lcm(M,p^fi-1));
\\COMPUTATION OF alpha^M:
alphaM=alpha;valphaM=min(valuation(lift(alphaM),p),0);
alphaM=Mod(Mod(lift(alphaM),mod*p^-valphaM),P)^M;
pi=Mod(lift(lift(alphaM))-1,P);if(pi==0,pi=Mod(mod,P));
\\COMPUTATION OF hatt:
hatt=1;
for(i=1,g,Primei=Ideal[i];Deltai=idealval(k,pi,Primei);rhoi=Lrho[i];
hatx=floor((log(rhoi)-log(Deltai*log(p)))/log(p))+1;
n=max(1,hatx);while(p^n*Deltai-rhoi*(N+n)<0,n=n+1);
hatn=n;if(hatn==1,hatti=N);if(hatn>1,hatti=p^hatn;
while(hatti*Deltai-(valuation(hatti,p)+N)*rhoi>0,hatti=hatti-1));
hatt=max(hatt,hatti));
\\COMPUTATION OF Log(alpha):
Log=pi;t=1;pit=pi;while(t<=hatt,t=t+1;pit=pi*pit;
Log=Log+(-1)^(t+1)*pit/t;
\\REDUCTION MODULO mod:
vLog=min(valuation(Log,p),0);Log=lift(Log)*p^-vLog;
Log =Mod(Log,mod*p^-vLog);Log=Mod(lift(Log)*p^vLog,P));
print("p=",p," hatt=",hatt);print("Log=",Log);for(j=1,g,
valLog=idealval(k,Log,Ideal[j]);print("valLogj=",valLog))}

p=2  hatt=128
Log=Mod(7*x^15+1/2*x^14+x^12+7*x^11+5*x^10+3*x^9 +63/8*x^8+x^7
   +6*x^5+19/4*x^4+2*x^3+9/2*x^2+9/2*x + 1/2, x^16-32*x^4+8*x-2)
valLogi=-40
time = 1,191 ms. ms.
\end{verbatim}\ns

The number ${\sf Log}$ is not a $2$-integer; the ramification index of 
${\mathfrak p} \mid 2$ is $16$ (too large, see Proposition \ref{prop1}\,(c)), 
whence $(\log(\alpha^M)) = \big[(4) \cdot {\mathfrak p}^8 \big]^{-1}$.

\smallskip
It is interesting to see the progression of the ${\mathfrak p}$-valuation
of ${\sf pi^t/t}$, as ${\sf t}$ increases, regarding ${\sf mod=2^3=Primei^{48}}$;
we have ${\sf pi = Mod(x^3+15*x+14,x^16-32*x^4+8*x-2)}$ of valuation $1$. 
This illustrates Proposition \ref{prop1} and Remark \ref{decrease} with another
case:

\smallskip
\ft\begin{verbatim}
-47,-62,-45,-76,-43,-58,-41,-88,-39,-54,-37,-68,-35,-50,-33,-96,-31,-46,-29,
-60,-27,-42,-25,-72,-23,-38,-21,-52,-19,-34,-17,-96,-15,-30,-13,-44,-11,-26,
-9,-56,-7,-22,-5,-36,-3,-18,-1,-6,1,-14,3,-28,5,-10,7,-40,9,-6,11,-20,13,-2,
15,-80,17,2,19,-12,21,6,23,-24,25,10,27,-4,29,14,31,-32,33,18,35,4,37,22,39,
-8,41,26,43,12,45,30,47,-32,49,34,51,20,53,38,55,8,57,42,59,28,61,46,63,0,
65,50,67,36,69,54,71,24,73,58,75,44,77,62,79,-32,81,66,83,52,85,70,87,4,...
\end{verbatim}\ns

\smallskip
The summation is significant, but many negative valuations will give the 
denominator $8$, and the necessary number $16$ appears at the end. 

\smallskip
The following example with a large prime $p$ is general and the running needs 
less than 3 seconds, noting that most of the time is spent by ${\sf bnfinit(P,1)}$:

\smallskip
\ft\begin{verbatim}
allocatemem(800000000)
{P=x^17-31*x-1;p=100000007;N=2;mod=p^N;k=bnfinit(P,1);
alpha=Mod(x,P);Ideal=List;Lrho=List;Lf=List;LDelta=List;LM=List;
F=idealfactor(k,p);g=matsize(F)[1];
for(j=1,g,Primej=component(F,1)[j];listput(Ideal,Primej);
rhoj=Primej[3];fj=Primej[4];listput(Lrho,rhoj);listput(Lf,fj);listput(LM,p^fj-1));
M=1;for(j=1,g,M=lcm(M,LM[j]));alphaM=alpha;
valphaM=min(valuation(lift(alphaM),p),0);
alphaM=Mod(Mod(lift(alphaM),mod*p^-valphaM),P)^M;
pi=Mod(lift(lift(alphaM-1)),P);if(pi==0,pi=Mod(mod,P));
for(j=1,g,Primej=Ideal[j];Deltaj=idealval(k,pi,Primej);listput(LDelta,Deltaj));
hatt=1;
for(j=1,g,rhoj=Lrho[j];hatx=floor((log(rhoj)-log(Deltaj*log(p)))/log(p))+1;
n=max(1,hatx);while(p^n*Deltaj-rhoj*(N+n)<0,n=n+1);
hatn=n;if(hatn==1,hattj=N);if(hatn>1,hattj=p^hatn;
while(hattj*Deltaj-(valuation(hattj,p)+N)*rhoj>0,hattj=hattj-1));
hatt=max(hatt,hattj));print("hatt=",hatt);
Log=pi;t=1;pit=pi;while(t<=hatt,t=t+1;pit=pi*pit;Log=Log+(-1)^(t+1)*pit/t;
vLog=min(valuation(lift(Log),p),0);Log=lift(Log)*p^-vLog;
Log=Mod(Log,mod*p^-vLog);Log=Mod(lift(Log)*p^vLog,P));
print("valuations(pi)=",LDelta);
print("ramifications=",Lrho);print("residue degrees=",Lf);
print("LM=",LM);print("M=",M);print("pi=",pi);print("Log=",Log)}

hatt=2  valuations(pi)=List([1,1,1,1,1])  ramifications=List([1,1,1,1,1])  
residue degrees=Lf=List([1,1,2,4,9])
LM=List([100000006,100000006,10000001400000048,100000028000002940000137200002400,
1000000630000176400028812003025260211768209882516296475485188320940353606])
M=100000085000033120007822601247358141463991700553911137815521697563769359
1518786629056793341442400
pi=Mod(3409200238644*x^16+5484595083921629*x^15+2119742948381996*x^14
+8987355229114822*x^13+7850508049535525*x^12+7732600941282028*x^11
+7477157623400997*x^10+2088556446198941*x^9+1916618434163281*x^8
+6251496237604706*x^7+400278228019474*x^6+326693722868559*x^5
+8452449591671430*x^4+6869203380844203*x^3+9922633494584296*x^2
+439241030746870*x+6959356087154892,x^17-31*x-1)
Log=3409200238644*x^16+5484595083921629*x^15+2119742948381996*x^14
+8987355229114822*x^13+7850508049535525*x^12+7732600941282028*x^11
+7477157623400997*x^10+2088556446198941*x^9+1916618434163281*x^8
+6251496237604706*x^7+400278228019474*x^6+326693722868559*x^5
+8452449591671430*x^4+6869203380844203*x^3+9922633494584296*x^2
+439241030746870*x+6959356087154892
time = 1,925 ms
\end{verbatim}\ns

\section{Programs testing \texorpdfstring{$\CH_k^\infty \ne 1$}{Lg}
for any number field}

The $\Log$-function allows to determine the structures of 
$\Gal(\wt k \cap H_k^\nr/k)$ and $\CH_k^\infty$; this leads to an elaborate 
program, unnecessary since our method is to test if $\CH_k^\infty = 1$ or not.
We first give the case of an imaginary quadratic field, then we
generalize to an arbitrary number field of any signature.

\begin{remark}
If, in an arbitrary number field $k$, the $p$-rank of $\CH_k$ is strictly larger 
than $r_2+1$, then $k$ is non-$p$-rational. For imaginary quadratic fields and 
$p>2$ it is well known that the part of the $p$-Hilbert class field contained 
in $\wt k$ is in fact contained in the anti-cyclotomic $\Z_p$-extension $k^\acyc$,
because of the semi-simple action of the complex conjugation $\tau$ on
$\Gamma := \Gal(\wt k/k)$, giving $\Gamma = \Gamma^+ \oplus \Gamma^-
\simeq \Gal(k^\cyc/k) \oplus \Gal(k^\acyc/k)$ in an obvious
manner (see \cite{Gra2026} for programs giving the first layer of 
$k^\acyc$), so $k$ is non-$p$-rational as soon as $\rk_p(\CH_k) \geq 2$. 
For instance, consider the first known imaginary quadratic fields with $\rk_p(\CH_k) 
\geq 3$ (Diaz y Diaz \cite[Appendice 1]{DYD1974}); for all them, $\CW_k^\bp = 1$
and Program \ref{P1} confirms the non-$3$-rationality and the fact that 
$\rk_3(\CH_k^\infty) = 2$. Up to $10^7$ only the case $-m = -5067967$ 
gives a non-$p$-rationality for a prime larger than $3$:

\smallskip
\ft\begin{verbatim}
x^2+5067967   
class group=[69,3,3]   Pbad=List([2,3,23,5067967])
**********************
p=2 Structure of the 2-torsion group:Tk=List([2])
k is not 2-rational
**********************
p=3 Structure of the 3-torsion group:Tk=List([3,3])
k is not 3-rational
**********************
p=23 Structure of the 23-torsion group:Tk=List([23])
k is not 23-rational
**********************
\end{verbatim}\ns
\end{remark}

\subsection{Case of an imaginary quadratic field}
The main characteristic is that there is no normalized $p$-adic regulator.

\subsubsection{Generalities}
We have the identity $\order \CT_k = \order \CH_k^\infty \cdot \order \CW_k^\bp$; 
the computation of $\order \CW_k^\bp$ is immediate, that of $\CH_k^\infty$ being 
much more complex since $\Gal(\wt k/k) \simeq \Z_p^2$, which needs the use of 
the $\Log$-function.
Of course, since the number of primes $p$ dividing $\order \CH_k$ is finite,
$\CT_k = 1$ for almost all $p$. This conclusion is only obvious for imaginary
quadratic fields; as soon as there are units of infinite order, this constitutes 
a conjecture for regulators of algebraic numbers, given in many of our papers.

\smallskip
For $p > 2$, one may assume that $\CH_k$ is cyclic; since $\BN_{k/\Q}
(\CH_k) = 1$, the action of the complex conjugation $\tau$ gives
$\CH_k = \CH_k^-$. Thus, $\order \CH_k^\infty = \ffrac{\order \CH_k}
{\big[\langle \frac{1}{n} \log(\alpha)^- \rangle + \log(\CU_k)^- : \log(\CU_k)^-\big]}$, 
where $\log(\alpha)^- = \frac{1-\tau}{2}\cdot \log(\alpha)$ is of the form
${\sf Mod(a \cdot x, x^2 + m)}$; whence, from ${\mathfrak a}^n = (\alpha)$:
$$\Log ({\mathfrak a})^- = \hbox{$\frac{1}{n}$} \log(\alpha)^-
 = a \cdot  \sqrt{-m}, \ \, a \in p^\Z \cdot \Z_p^\times .$$

So, $\log(\CU_k^-)$ being of the form $b \cdot \sqrt{-m} $, $b \in \Z_p$, 
the $p$-valuation of $\order \CH_k^\infty$ is given in the program by
$v_p(\order \CH_k) - v_p([\wt k \cap H_k^\nr : k])$, thus by:
$${\sf valuation(hk,p) - \big[valuation(b,p) - \min(valuation(a,p),
valuation(b,p)\big]}. $$

The case $p = 2$ is particular since the action of $\tau$ is not semi-simple.
So the test of the $2$-rationality of $k = \Q(\sqrt{-m})$ with the ray class 
group of modulus $8$, is faster.

\subsubsection{Program computing \texorpdfstring{$\CH_k^\infty$ 
for $k = \Q(\sqrt{-m})$}{Lg}} \label{quadHinfty}

For $p \ne 2$, we know that $\log(\CU_k)^- = p \cdot \Z_p \sqrt{-m}$ if $p$
is unramified in $k$ (i.e. $p \nmid m$), and $\log(\CU_k)^- =  \Z_p  \sqrt{-m}$
if $p \mid m$, $-m \not \equiv -3 \pmod 9$ if $p = 3$, then
$\log(\CU_k)^- = 3 \cdot \Z_3 \sqrt{-m}$ if $p = 3$ and $-m \equiv -3 \pmod 9$
(cf. \cite[Lemma 2.6]{Gra2026}). This defines the number $b$.

\smallskip
In the {\sc pari/gp} notation, a canonical generator of the $\Z_p$-module
$\log(\CU_k)^-$ is ${\sf Mod(p*x,P)}$ if $p \nmid m$, and ${\sf Mod(x,P)}$
if $p \mid m$.

The program does not consider, for $p = 3$, the (numerous) cases where 
$\CW_k^\bp(3) = 3$.
A sufficient precision is necessary, only for few values of $m$, because of the 
instruction ${\sf bnfisprincipal(k,C0)[2]}$; recall that here, ${\sf mod = p^{n0+2}}$
to take into account the relation $\log({\mathfrak a}) = \frac{1}{n_0}\log(\alpha_0)$,
modulo $p^2$:

\smallskip
\ft\begin{verbatim}
allocatemem(800000000)
\p 550
{minm=2;maxm=10^4;minp=3;maxp=10^4;
for(m=minm,maxm,if(core(m)!=m,next);P=x^2+m;k=bnfinit(P);
\\TEST ON W:
if(-m!=-3 & Mod(-m,9)==-3,next);
\\TESTS ON THE CYCLICITY OF THE p-CLASS GROUP:
Hk=k.clgp;Nk=Hk[2];A=Hk[3];dA=matsize(A)[2];
forprime(p=minp,maxp,valhkp=valuation(k.no,p);if(valhkp==0,next);
Hkp=List;dp=0;for(i=1,dA,vNk=valuation(Nk[i],p);
if(vNk!=0,dp=dp+1;listput(Hkp,p^vNk)));if(dp>1,next);
\\THE p-CLASS GROUP IS NOW p-CYCLIC (A0,N0):
for(i=1,dA,if(Mod(Nk[i],p)==0,N0=Nk[i];A0=A[i];break));
n0=valuation(N0,p);N=n0+2;mod=p^N;
\\DECOMPOSITION OF p:
F=idealfactor(k,p);g=matsize(F)[1];Prime=component(F,1)[1];
\\PRIME-TO-p GENERATOR B0 OF THE p-CLASS GROUP:
B0=A0;Norm0=idealnorm(k,B0);w0=valuation(Norm0,p);if(w0!=0,
\\RAMIFIED CASE:
if(Mod(m,p)==0,z=idealappr(k,B0);B0=idealdiv(k,z,B0));
\\UNRAMIFIED CASE:
if(Mod(m,p)!=0,Bp=idealmul(k,B0,p);z=idealappr(k,Bp);
B0=idealdiv(k,z,Bp)));
\\COMPUTATION of THE VALUATION OF Log(U)^-: 
if(p!=3,if(Mod(m,p)!=0,vLogU=1);if(Mod(m,p)==0,vLogU=0));
if(p==3,if(Mod(m,3)!=0,vLogU=1);if(Mod(m,3)==0,vLogU=0));
\\DECOMPOSITION OF (p)--COMPUTATION OF M (p-1 or p^2-1):
rho=Prime[3];f=Prime[4];if(f==1,M=p-1);if(f==2,M=p^2-1);
\\COMPUTATION OF beta0 AND alpha0M:
C0=idealpow(k,B0,N0);beta0=bnfisprincipal(k,C0)[2];
alpha0=Mod(beta0[1]*k.zk[1]+beta0[2]*k.zk[2],P);
\\REDUCTION OF alpha0^M MODULO mod:
alpha0M=Mod(Mod(lift(alpha0),mod),P)^M;
pi=Mod(lift(lift(alpha0M))-1,P);if(pi==0,pi=Mod(mod,P));
\\COMPUTATION OF hatt AND THE LOG SERIES MODULO mod:
hatt=1;for(j=1,g,Primej=component(F,1)[j];
Deltaj=idealval(k,pi,Primej);rhoj=Primej[3];
hatx=floor((log(rhoj)-log(Deltaj*log(p)))/log(p))+1;
n=max(1,hatx);while(p^n*Deltaj-rhoj*(N+n)<0,n=n+1);
hatn=n;if(hatn==1,hattj=N);if(hatn>1,hattj=p^hatn;
while(hattj*Deltaj-(valuation(hattj,p)+N)*rhoj>0,hattj=hattj-1));
hatt=max(hatt,hattj));
Log0=pi;t=1;pit=pi;Log0=pi;t=1;pit=pi;
\\COMPUTATION OF Log0=Log(alpha0M):
while(t<=hatt,t=t+1;pit=pi*pit;Log0=Log0+(-1)^(t+1)*pit/t);
Log0=Mod(lift(Log0),mod);a=polcoeff(Log0,1);
if(a==0,vala=N+2);if(a!=0,vala=valuation(a,p));
vlog=vala-n0;Minval=min(vLogU,vlog);
Index=vLogU-Minval;valhkinfty=valhkp-Index;if(valhkinfty>0,
print("-m=",-m," k non ",p,"-rational"," Index=",p^Index,
" Hkp=",Hkp," hkinfty=",p^valhkinfty))))}
time = 25,816 ms.
(...)
-m=-9527 k non 11-rational Index=1 Hkp=List([11]) hkinfty=11
-m=-9529 k non 7-rational  Index=1 Hkp=List([7])  hkinfty=7
-m=-9553 k non 5-rational  Index=1 Hkp=List([5])  hkinfty=5
-m=-9554 k non 5-rational  Index=1 Hkp=List([5])  hkinfty=5
-m=-9557 k non 13-rational Index=1 Hkp=List([13]) hkinfty=13
-m=-9587 k non 23-rational Index=1 Hkp=List([23]) hkinfty=23
-m=-9593 k non 5-rational  Index=1 Hkp=List([5])  hkinfty=5
-m=-9622 k non 3-rational  Index=1 Hkp=List([9])  hkinfty=9
-m=-9626 k non 3-rational  Index=3 Hkp=List([9])  hkinfty=3
-m=-9634 k non 7-rational  Index=1 Hkp=List([7])  hkinfty=7
-m=-9663 k non 3-rational  Index=9 Hkp=List([27]) hkinfty=3
\end{verbatim}\ns

\subsection{Case of an arbitrary number field}
In that case, the three $\Z_p$-modules $\CH_k^\infty$, $\CR_k$ and 
$\CW_k^\bp$ may be non-trivial; we refer to  Programs \ref{P5} and 
\ref{PW} for independent tests of $\CR_k$ and $\CW_k^\bp$, respectively.

\subsubsection{Generalities and hypothesis}
We assume that we know that $\CW_k^\bp = 1$ (otherwise $k$ is non
$p$-rational) and that $\CH_k \ne 1$ (otherwise $\CH_k^\infty = 1$).

\smallskip
We will use the method of determining the small bad primes 
(including $p = 2$) to which we apply
the usual classical computation of $\rk_p(\CT_p)$; this simplifies the $\log$  
and matrix computations. Then we must first compute $\CR_k = 
\tor_{\Z_p} (\CU_k/\ov E_k) \simeq \tor_{\Z_p}
\big (\log(1 + p \cdot \CO_k) \big /\log(\ov E_k) \big)$.

\smallskip
Once this is done, with $\CR_k$ invertible and $\CW_k = 1$, it remains 
the most difficult invariant $\CH_k^\infty$ where some simplifications 
arise (non-ramification of $p$, $\CR_k = 1$ meaning that $\log(\ov E_k)$, 
is direct factor in $\log(\CU) = p \cdot \CO_k$ since $p \ne 2$). This point is 
essential since the image of $\Log$ takes place in the $\Q_p$-vector 
space $\oplus_{v \mid p} k_v /\Q_p \cdot \log(\ov E_k)$.

\smallskip
It is clear that $\{1, x, \ldots , x^{d-1}\}$ is a $\Z_p$-basis of $\CO_k$ 
when $p$ is not bad, especially if $p$ does not divide the discriminant 
$\CD_k$, nor the denominators of integer basis elements; by security, 
we may compute the corresponding matrix determinant; indeed, it is not 
certain that primes dividing the polynomial discriminant are non-bad. 
The integer basis of $k$ is given by ${\sf k.zk}$, where
${\sf k.zk[i] = Mod(c_0^i + c_1^i*x + \ldots, c_{d-1}^i*x^{d-1},P)}, \ 
\hbox{${\sf i = 1,\ldots , d}$}$. 
So we compute the $p$-valuation of $\det(c_j^i)$:

\smallskip
\ft\begin{verbatim}
{P=x^9+7*x^3-1400*x-7;maxp=10^8;d=poldegree(P);if(polisirreducible(P)==0,break);
k=bnfinit(P);Dk=k.disc;print(factor(Dk));denB=1;for(j=1,d,Bzk=lift(k.zk[j]);
for(i=0,d-1,denB=lcm(denominator(polcoeff(Bzk,i)),denB)));
forprime(p=2,maxp,if(Mod(Dk*denB,p)==0,print("bad prime p=",p);next);
Mz=List;for(j=1,d,zj= k.zk[j];Lc=List;for(i=0,d-1,coeff=Mod(polcoeff(zj,i),p);
listput(Lc,coeff));listput(Mz,Lc));MZ=matrix(d,d,u,v,Mz[u][v]);
valdet=valuation(matdet(MZ),p);if(valdet!=0,
print("p=",p," warning"," valdeterminant=",valdet)))}

Dk=5^3*7^8*249517*6670822160235662179
bad prime p=5
bad prime p=7
bad prime p=17
bad prime p=249517
time = 54,537 ms for p<10^7.
time = 7min, 52,584 ms for p<10^8.
\end{verbatim}\ns

The prime $p = 17$ is not ramified but divides the polynomial discriminant,
and the denominator of several elements of ${\sf k.zk}$, as:
\ft$${\sf \frac{3}{17}*x^8-\frac{2}{17}*x^7+\frac{7}{17}*x^6
+\frac{1}{17}*x^5+\frac{5}{17}*x^4+\frac{8}{17}*x^3
+\frac{27}{17}*x^2-\frac{18}{17}*x-\frac{3695}{17}}. $$\ns

\subsubsection{Test for \texorpdfstring{$\CH_k^\infty \ne 1$}{Lg} 
for any number field}

We assume that the normalized $p$-adic regulator has been computed via 
Program \ref{P5} and is invertible (otherwise, $k$ is non-$p$-rational).

\smallskip
Because of the instruction ${\sf betai=bnfisprincipal(k,Bi)[2]}$, giving the 
components on the integer basis of the pseudo-units $\alpha_i$ generating 
${\mathfrak a}_i^{n_i}$, it is crucial to take a sufficient precision. 
The program gives the $d$ components of $\Log({\mathfrak a}_i)$ 
modulo $\Log(P_k \otimes \Z_p) = p\cdot \CO_k + \Q_p \cdot \log(\ov E_k)$; 
then, we have the following lemma:

\begin{lemma}\label{bimatrix}
Assume that $p > 2$ is unramified in $k$, does not divide the denominators
of the integer basis elements of $k$, and that $\CR_k = 1$. 
Put $(\alpha_i) = {\mathfrak a}_i^{n_i}$, $n_i \equiv 0\! \pmod p$, where the 
$\Ccl({\mathfrak a}_i)$ are the generating classes given by {\sc pari/gp} in 
${\sf k.clgp}$.
Let $\{\alpha_1, \ldots, \alpha_h; \varepsilon_1, \ldots, \varepsilon_r\}$,
$h = \rk_p(\CH_k)$, $r := \rk_\Z(E_k)$, 
be a fundamental system of pseudo-units of $k$. Set:
\begin{equation*}
\left \{ \begin{aligned}
\ffrac{1}{p} \log(\alpha_i) & = \sm_{\ell=0}^{d-1} a_\ell^i \cdot x^\ell 
\in \Z_p[x],\ \,  i \in [1, h], \\\
\ffrac{1}{p} \log(\varepsilon_j) & = \sm_{\ell=0}^{d-1} b_\ell^i \cdot x^\ell 
 \in \Z_p[x],\ \, j \in [1, r].
\end{aligned}\right.
\end{equation*}

Let $\CX_k$ be the $(d, h+r)$-matrix defined by the 
coefficients $a_\ell^i$ and $b_\ell^j$ modulo $p$.
Then, the invariant $\CH_k^\infty$ is non-trivial if and only if the rank 
of the $\F_p$-matrix $\CX_k$ is strictly less than $\rk_p(\CR_k) + \rk_p(\CH_k)
= \rk_\Z(E_k) + \rk_p(\CH_k) = r+h$ (${\sf rE + rhk}$ in the program).
\end{lemma}

\begin{proof}
Since $p > 2$ is unramified, we have $\log(\CU_k) = p \cdot \CO_k$.
We have $\CH_k^\infty \ne 1$ if and only if there exists a class 
$\Ccl({\mathfrak a})$ of order $p$ in $\CH_k$, ${\mathfrak a} \in I_k \otimes \Z_p$,
such that $\Log({\mathfrak a}) = 0$. Put:
\begin{equation} \label{Ap=1}
{\mathfrak a} = \prd_{i=1}^h {\mathfrak a}_i^{\lambda_i \cdot n_i/p},
\end{equation} 

\noindent
up to a principal ideal in $P_k \otimes \Z_p$, with $\lambda_i \in \F_p$ 
not all zero (this means that the class of  ${\mathfrak a}$ is of order~$p$). 
Since $\log(\CU_k) = p \cdot \CO_k$ and by assumption:
\begin{equation} \label{tor=0}
\CR_k \simeq \tor_{\Z_p}(\log(\CU_k)/\log(\ov E_k)) = 1, 
\end{equation} 

\noindent
$\log(\ov E_k)$ is a direct factor in $p \cdot \CO_k$; the relation 
$\Log({\mathfrak a})  = 0$ is equivalent to $\Log({\mathfrak a} \cdot
(u)^{-1}) = 0$ for some $(u) \in P_k \otimes \Z_p$, whence
$\Log({\mathfrak a} \cdot (u)^{-1}) \in \Q_p \cdot \log(\ov E_k)$ 
(Formula \eqref{alpha}); since ${\mathfrak a} \cdot (u)^{-1}$ 
(still denoted ${\mathfrak a}$ by abuse) represents the same 
class of order $p$, one may write:
\begin{equation} \label{z}
\sm_{i=1}^h \lambda_i \ffrac{1}{p} \log(\alpha_i) = p^z \cdot 
\ffrac{1}{p}\log(\varepsilon),\ \, z \in \Z,\  \varepsilon \in \ov E_k.
\end{equation} 

The left member is in $\CO_k$, so $p^z \cdot \ffrac{1}{p}
\log(\varepsilon) \in \CO_k$. Since, from the relation $\eqref{tor=0}$,
$\log(\varepsilon) \equiv 0 \pmod {p^{m+1}}$, $m\geq 0$, implies 
$\varepsilon = \eta^{p^m}$, $\eta \in \ov E_k$, on may suppose that 
$\ffrac{1}{p} \log(\varepsilon) \in \CO_k \setminus p \CO_k$; writing this 
on a basis of the free $\Z_p$-module $\CO_k$, we deduce that $z = 0$.
Whence, from \eqref{z}, the existence of elements $\lambda_i \in \F_p$, 
not all zero, and $\mu_j \in \F_p$, such that:
$$\sm_{i=1}^h \lambda_i \ffrac{1}{p} \log(\alpha_i) + 
\sm_{j=1}^r \mu_j \ffrac{1}{p} \log(\varepsilon_j) \equiv 0 \pmod p, $$

\smallskip
\noindent
so the $\F_p$-rank of the matrix
$\CX_k$ is strictly less that $\rk_\Z(E_k) + \rk_p(\CH_k) = r+h$.
Since $\CR_k$ is invertible by \eqref{tor=0}, the rank of $\CX_k$ is
at least $r$ (the program verifies that $\CR_k$ is invertible; if not
it concludes that $k$ is not $p$-rational for another reason).
\end{proof}

\subsubsection{Program testing $\CH_k^\infty$ for any
number field and non bad primes}\label{P6}
The hypothesis are those given in Lemma \ref{bimatrix}; in particular,
that it has been verified that $\CR_k = 1$.

\smallskip
\ft\begin{verbatim}
allocatemem(8000000000)
\p 1000
{for(a=100,200, P=x^8-7*a*x^4-1;bp=3;Bp=10^4;
if(polisirreducible(P)==0,next);d=poldegree(P);k=bnfinit(P,1);
kP=nfinit(P);r1=kP[2][1];r2=kP[2][2];rE=r1+r2-1;rk=r2+1;
\\BAD PRIMES (ramified + divisors of denominators of k.zk):
Pbad=List;
\\DENOMINATORS OF THE INTEGER BASIS and Dk:
denB=1;for(j=1,d,Bzk=lift(k.zk[j]);for(i=0,d-1,
denB=lcm(denominator(polcoeff(Bzk,i)),denB)));
Fbad=factor(abs(k.disc)*denB);dF=matsize(Fbad)[1];
Prodbad=1;for(i=1,dF,p=component(Fbad,1)[i];
listput(Pbad,p);Prodbad=Prodbad*p);
\\TEST FOR THE NON-BAD PRIMES:
\\DECOMPOSITION OF p IN k, COMPUTATION OF M:
forprime(p=3,Bp,if(Mod(Prodbad,p)==0,next);
valhkp=valuation(k.no,p);if(valhkp==0,next);print();
print("a=",a," P=",P);print("class group=",k.cyc);print("Pbad=",Pbad);
print("number of independent Z_p-extensions=",rk);
print("Z-rank of the group of units=",rE);
\\ p-RANK OF Hk, LIST OF PRIME-TO-p IDEALS:
Hk=k.clgp;Nk=Hk[2];A=Hk[3];dA=matsize(A)[2];B=List;
N=0;rhk=0;for(i=1,dA,Nui=Nk[i];vi=valuation(Nui,p);
if(vi>0,rhk=rhk+1;N=max(N,vi)+2;CA=A[i];nu=idealnorm(k,CA);
if(Mod(nu,p)==0,Y=idealmul(k,CA,p);z=idealappr(k,Y);
Z=idealdiv(k,z,Y);listput(B,Z));if(Mod(nu,p)!=0,listput(B,CA))));
mod=p^(N+3);F=idealfactor(k,p);g=matsize(F)[1];M=1;
for(i=1,g,fi=component(F,1)[i][4];M=lcm(M,p^fi-1));
\\COMPUTATION OF THE LOG(alphai):
LogA=List;for(i=1,rhk,Bi=idealpow(k,B[i],Nk[i]);
betai=bnfisprincipal(k,Bi)[2];alphai=0;
for(j=1,d,alphai=alphai+Mod(betai[j]*k.zk[j],P));
\\REDUCTION OF alphai^M MODULO mod:
alphaiM=Mod(Mod(lift(alphai),mod),P)^M;
pi=Mod(lift(lift(alphaiM)),P)-1;if(pi==0,pi=Mod(mod,P));
\\COMPUTATION OF hatt FOR Log(alphaiM):
hatt=1;
for(j=1,g,Primej=component(F,1)[j];rhoj=Primej[3];
Deltaj=idealval(k,pi,Primej);rhoj=Primej[3];
hatx=floor((log(rhoj)-log(Deltaj*log(p)))/log(p))+1;
n=max(1,hatx);while(p^n*Deltaj-rhoj*(N+n)<0,n=n+1);
hatn=n;if(hatn==1,hattj=N);if(hatn>1,hattj=p^hatn;
while(hattj*Deltaj-(valuation(hattj,p)+N)*rhoj>0,hattj=hattj-1));
hatt=max(hatt,hattj));
Log=pi;t=1;pit=pi;
while(t<=hatt,t=t+1;pit=pi*pit;Log=Log+(-1)^(t+1)*pit/t);
listput(LogA,lift(Log)));
\\COMPUTATION OF THE Log(ai)=1/ni*log(alphai)
MA=List;for(ell=1,rhk,Lell=LogA[ell]/p;
Lc=List;for(i=0,d-1,coeff=polcoeff(Lell,i);
listput(Lc,coeff));listput(MA,Lc));
MLogA=matrix(rhk,d,u,v,Mod(MA[u][v],p));
\\TEST ON THE REGULATOR:
E=k.fu;LogE=List;for(ell=1,rE,
epsilonellM=Mod(Mod(lift(E[ell]),mod),P)^M;
pi=Mod(lift(lift(epsilonellM)),P)-1;
\\COMPUTATION OF hatt FOR Log(epsilonellM):
N0=2;if(p==2,N0=3);
hatt=1;for(j=1,g,Primej=component(F,1)[j];rhoj=Primej[3];
Deltaj=idealval(k,pi,Primej);rhoj=Primej[3];
hatx=floor((log(rhoj)-log(Deltaj*log(p)))/log(p))+1;
n=max(1,hatx);while(p^n*Deltaj-rhoj*(N0+n)<0,n=n+1);
hatn=n;if(hatn==1,hattj=N);if(hatn>1,hattj=p^hatn;
while(hattj*Deltaj-(valuation(hattj,p)+N)*rhoj>0,hattj=hattj-1));
hatt=max(hatt,hattj));
Log=pi;t=1;pit=pi;
while(t<=hatt,t=t+1;pit=pi*pit;Log=Log+(-1)^(t+1)*pit/t);
listput(LogE,lift(Log)));
ME=List;for(ell=1,rE,Lell=LogE[ell]/p;
Lc=List;for(i=0,d-1,coeff=Mod(polcoeff(Lell,i),p);
listput(Lc,coeff));listput(ME,Lc));
MLogE=matrix(rE,d,u,v,ME[u][v]);rankE=matrank(MLogE);
if(rankE<rE,print(p,"-regulator non invertible,",
" k is not ",p,"-rational: Hinfty>1");next);
X=concat(MLogA~,MLogE~);rankX=matrank(X);
print("p=",p," rank(R)=",matrank(MLogE)," rank(X)=",matrank(X));
if(rankX==rE+rhk,print("Hinfty=1"));
if(rankE<rE || rankX<rE+rhk,
print("k is not ",p,"-rational: Hinfty>1"))))}
time = 1min, 35,921 ms.
\end{verbatim}\ns

Let's give some examples:

\subsubsection{Case $a = 5$}
Program \ref{P5} shows that $\CR_k = 1$:

\smallskip
\ft\begin{verbatim}
a=5 P=x^8 - 35*x^4 - 1  
class group=[3,3]
Pbad=List([2,1229])
number of independent Z_p-extensions=4
Z-rank of the group of units=4
p=3 rank(R)=4 rank(X)=5
k is not 3-rational: Hinfty>1
\end{verbatim}\ns

\smallskip
Since $3$ is not a bad prime, we can describe logarithms modulo $3$ 
on the $3$-basis $\{1, x,\ldots, x^7\}$ to find the various $\F_3$-relations, 
knowing $\log{\ov E_k}$ (direct factor of $\Z_3$-dimension $4$, of $\CO_k$
of $\Z_3$-dimension $8$). The matrix $\CX_k$ has the following form, 
where, for short, ${\sf logai = \frac{1}{p} \cdot \log(\alpha_i^M)}$ and 
${\sf logej = \frac{1}{p} \cdot \log(\varepsilon_j^M)}$, seen modulo $p$
on the $\Z_p$-basis $\{1, x, \ldots, x^{d-1}\}$:

\smallskip
\ft\begin{verbatim}
pseudo-units:
loga1=[1,1,1,0,2,0,0,2]
loga2=[0,2,2,0,0,0,0,1]
units:
loge1=[2,0,0,0,1,0,0,0]
loge2=[0,1,0,0,0,2,0,1]
loge3=[1,2,1,1,2,2,0,0]
loge4=[1,1,1,2,2,1,0,0]
\end{verbatim}\ns

\smallskip
One verifies easily that ${\sf loga1+loga2 = [1, 0, 0, 0, 2, 0, 0, 0] = 2*loge1}$
which means that, modulo $\Q_3 \log(E_k)$, the image of the $3$-class 
of ${\mathfrak a}_1 \cdot {\mathfrak a}_2$ is trivial (whence $H_k^\nr$ is not 
linearly disjoint from the compositum of the $\Z_3$-extensions of $k$). 

\smallskip
To see that $H_k^\nr$ is not contained in $\wt k$, we note for instance that 
$\Log({\mathfrak a}_1)$ is of order $3$ in $\Log(I_k \otimes \Z_p)/\Log
(P_k \otimes \Z_p)$. A direct computation of $\CT_k$ confirms the above
calculations since, from Program \ref{P2} with ${\sf nup > 2}$,
${\sf Tk=List([3])}$; whence $\CT_k \simeq \Z/3\Z$.

\subsubsection{Case $a = 90$}

\smallskip
\ft\begin{verbatim}
a=90 P=x^8 - 630*x^4 - 1  
class group=[94]
Pbad=List([2, 49613])
number of independent Z_p-extensions=4
Z-rank of the group of units=4
p=47 rank(R)=4 rank(X)=4
k is not 47-rational: Hinfty>1

pseudo-unit:
loga=[5,19,36,14,44,19,5,46]
units:
loge1=[33,0,0,0,46,0,0,0]
loge2=[0,16,0,21,0,0,0,17]
loge3=[37,8,4,34,6,0,11,32]
loge4=[0,35,0,44,0,27,0,38]
\end{verbatim}\ns

\smallskip
One finds ${\sf loga = 10*loge1+23*loge2+9*loge3+39*loge4}$ in $\F_{47}$,
proving that the image of the $p$-class group in $\Gal(\wt k/k)$
is not trivial; so that $k$ is not $47$-rational.

\begin{remark}
It is possible that some $\ffrac{1}{p} \log(\alpha)$ be zero modulo $p$ (whence 
non-nullity of the ${\mathfrak p}$-Fermat quotients $\delta_{\mathfrak p}(\alpha)$);
this happens rarely. We have the following example with another 
polynomial and $p = 3$:

\smallskip
\ft\begin{verbatim}
P=x^4-438*x^3+4*x-1  
class group=[6,2,2]
Pbad=List([2,5,29,563,148861])
number of independent Z_p-extensions=2
Z-rank of the group of units=2
p=3  rank(R)=2  rank(X)=2
k is not 3-rational: Hinfty>1
\end{verbatim}\ns

\smallskip
\noindent
Indeed, verifying independently, one finds that the ideal class 
$\Ccl({\mathfrak a})$, of order $6$, gives an element $\alpha$
of the form (after reduction modulo $81$):

\smallskip
\ft\begin{verbatim}
alpha=Mod(Mod(37,81)*x^3+Mod(75,81)*x^2+Mod(30,81)*x+Mod(37,81),P)
Since M=80, one finds:
alpha^M-1=Mod(Mod(27,81)*x^3+Mod(9,81)*x^2+Mod(54,81)*x+Mod(27,81),x^4-438*x^3+4*x-1)

pseudounit:
loga=[Mod(0,3) Mod(0,3) Mod(0,3) Mod(0,3)]
units:
loge1=[Mod(0,3) Mod(2,3) Mod(1,3) Mod(0,3)]
loge2=[Mod(0,3) Mod(0,3) Mod(0,3) Mod(1,3)]
\end{verbatim}\ns
\end{remark}

\subsubsection{Case $a = 21$}
In this case, ${\sf Hinfty>1}$ for various reasons:

\smallskip
\ft\begin{verbatim}
a=21 P=x^8-147*x^4-1 
class group=[273]
Pbad=List([2,21613])
number of independent Z_p-extensions=4
Z-rank of the group of units=4
3-regulator non invertible, k is not 3-rational:  Hinfty>1
7-regulator non invertible, k is not 7-rational:  Hinfty>1
p=13 rank(R)=4 rank(X)=4    k is not 13-rational: Hinfty>1
\end{verbatim}\ns

\smallskip
The case $a = 126$, $P=x^8 - 882*x^4 - 1$, needs precision ${\sf p 5000}$, 
due probably to the class number $1680 = 2^4 \cdot 3 \cdot 5 \cdot 7$, and 
then obtains that $k$ is not $p$-rational for $p = 3, 5, 7$ with ${\sf Hinfty>1}$.
Using Program \ref{P1}, one finds that $k$ is not $p$-rational for 
$p = 2, 3, 5, 7, 29, 281,$

\end{document}